\documentclass[12pt]{article}
\usepackage{amsmath,amssymb,amsfonts}
\usepackage{graphicx}
\usepackage{epsfig}
\usepackage{amsfonts}
\usepackage{amssymb, latexsym, amsmath, pb-diagram}
\usepackage{pstricks-add}
\usepackage{multicol}
\usepackage{bm}
\usepackage{mathrsfs}
\usepackage{xy}
\usepackage{epic,eepic}
\xyoption{all}
\begin{document}
\title{The Complement of Polyhedral Product Spaces and the Dual Simplicial Complexes}
\author {
Qibing Zheng\\School of Mathematical Science and LPMC, Nankai University\\
Tianjin 300071, China\\
zhengqb@nankai.edu.cn\footnote{Project Supported by Natural Science
Foundation of China, grant No. 11071125 and No. 11671154\newline
\hspace*{5.5mm}Key words and phrases: complement space; dual complex;
polyhedral product complex; universal algebra.\newline
\hspace*{5.5mm}Mathematics subject classification: 55N10}}\maketitle
\input amssym.def
\newsymbol\leqslant 1336
\newsymbol\geqslant 133E
\baselineskip=20pt
\def\w{\widetilde}
\def\RR{\mathscr R}
\def\SS{\mathscr S}
\def\XX{\mathscr X}

\begin{abstract} In this paper, we define and prove basic properties
of complement polyhedral product spaces, dual complexes and polyhedral product complexes.
Then we compute the universal algebra of polyhedral product complexes under certain split conditions
and the Alexander duality isomorphism on certain polyhedral product spaces.
\end{abstract}

\hspace*{40mm}${\displaystyle{\bf Table\,\,of\,\,Contents}}$

Section 1\, Introduction

Section 2\, Complement Spaces, Dual Complexes and Polyhedral Product Complexes

Section 3\, Homology and Cohomology Group

Section 4\, Universal Algebra

Section 5\, Duality Isomorphism
\vspace{3mm}

\newtheorem{Definition}{Definition}[section]
\newtheorem{Theorem}{Theorem}[section]
\newtheorem{Lemma}{Lemma}[section]
\newtheorem{Example}{Example}[section]
\font\hua=eusm10 scaled\magstephalf
\def\MM{\mathscr M}
\def\LL{\mathscr L}
\def\zz{\Bbb Z}

\section{Introduction}\vspace{3mm}

The polyhedral product theory, especially the homotopy type of polyhedral product
spaces, is developing rapidly nowadays. The first known polyhedral product space was the moment-angle
complex introduced by Buchstaber and Panov \cite{BP} and
was widely studied by mathematicians in the area of toric topology and geometry
(see \cite{A},\cite{B},\cite{GM},\cite{GR},\cite{GJ},\cite{GS}).
Later on, the homotopy types of polyhedral product spaces were studied by
Grbi\'{c} and Theriault \cite {GR},\cite{GJ},\cite{GS}, Beben and Grbi\'{c} \cite{GB},
Bahri, Bendersky, Cohen and Gitler \cite{B1},\cite{B2},\cite{B3} and many others
(\cite {BS},\cite{DJ},\cite{DE}).
The cohomology ring of homology split polyhedral product spaces and the
cohomology algebra over a field of polyhedral product spaces were computed in \cite{Z}.

In this paper, a polyhedral product space ${\cal Z}(K;\underline{X},\underline{A})$ with
$(\underline{X},\underline{A})=\{(X_k,A_k)\}_{k=1}^m$ is general than usual in that each $(X_k,A_k)$ is
a topological pair but not a CW-complex pair (see Definition 2.1).
Then for $M={\cal Z}(K;\underline{X},\underline{A})$,
is the complement space $M^c=(X_1{\times}{\cdots}{\times}X_m)\setminus M$ a polyhedral product space?
In Theorem 2.4, we show that $M^c={\cal Z}(K^\circ;\underline{X},\underline{A}^c)$, where $K^\circ$
is the dual complex of $K$ relative to $[m]$ and $A_k^c=X_k{\setminus}A_k$ is the complement space of $A_k$.

Let ${\cal Z}(K;\underline{Y},\underline{B})$,
$(\underline{Y},\underline{B})=\{(Y_k,B_k)\}_{k=1}^m$ be the polyhedral product space
defined as follows. For each $k$, $(Y_k,B_k)$ is a pair of polyhedral product spaces given by ($s_k=n_1{+}{\cdots}{+}n_k$)
$$Y_k={\cal Z}(X_k;\underline{U_k},\underline{C_k}),\,\, B_k={\cal Z}(A_k;\underline{U_k},\underline{C_k}),\,\,
(\underline{U_k},\underline{C_k})=\{(U_i,C_i)\}_{i=s_{k-1}{+}1}^{s_k},$$
where $(X_k,A_k)$ is a simplicial pair on $[n_k]$.
In theorem 2.9, we prove that
${\cal Z}(K;\underline{Y},\underline{B})$ is also a polyhedral product space
${\cal Z}({\cal S}(K;\underline{X},\underline{A});\underline{U},\underline{C})$,
where ${\cal S}(K;\underline{X},\underline{A})$ is the polyhedral product complex
defined in Definition 2.7.
When all $(X_k,A_k)=(\Delta\!^{n_k},K_k)$, the simplicial complex ${\cal S}(K;\underline{X},\underline{A})$ is
denoted by ${\cal S}(K;K_1,{\cdots},K_m)$ which is just
the composition complex $K(K_1,{\cdots},K_m)$ defined by Anton in Definition 4.5 of \cite{AY}.

In Section 3,
we compute the reduced (co)homology group and the (right) total (co)homology group of polyhedral product complexes
from the point of view of split inclusion (Theorem 3.9 and Theorem 3.11) .
In Example~3.10, we show that if $K$ and all $L_k$ are homology spheres, then
${\cal S}(K;L_1,{\cdots},L_m)$ is a homology sphere. This result is in accordance
with that of \cite{AY}, where the homotopy type of ${\cal S}(K;L_1,{\cdots},L_m)$
is studied by Anton.

In section~4, we compute cohomology algebra of a wide class of complexes in Theorem 4.5
including the cohomology algebra of ${\cal Z}(K;\underline{Y},\underline{B})$ mentioned above in Example 4.6.
In Theorem 5.6, we compute the Alexander duality isomorphism on the pair $
(X_1{\times}{\cdots}{\times}X_m,{\cal Z}(K;\underline{X},\underline{A}))$,
where all $X_k$'s are orientable manifolds and all $A_k$'s are polyhedra.
\vspace{3mm}

\section{Complement Spaces, Dual Complexes
and Polyhedral Product complexes}\vspace{3mm}

{\bf Conventions and Notations} For a finite set $S$, $\Delta\!^S$ is the simplicial complex with only one maximal simplex $S$,
i.e., it is the set of all subsets of $S$ including the empty set $\emptyset$.
Define $\partial \Delta\!^S=\Delta\!^S{\setminus}\{S\}$.
For $[m]=\{1,{\cdots},m\}$, $\Delta\!^{[m]}$ is simply denoted by $\Delta\!^m$.
Specifically, define
$\Delta\!^0=\Delta\!^\emptyset=\{\emptyset\}$ and $\partial\Delta\!^0=\{\,\}$.
The void complex $\{\,\}$ with no simplex is inevitable in this paper.

For a simplicial complex $K$ on $[m]$ (ghost vertex $\{i\}\notin K$ is allowed) and $\sigma\subset[m]$ ($\sigma\notin K$ is allowed),
the link of $\sigma$ with respect to $K$ is the simplicial complex
${\rm link}_K\sigma=\{\tau\,|\,\sigma{\cup}\tau\in K,\,\sigma{\cap}\tau=\emptyset\}$.
This implies ${\rm link}_K\sigma=\{\emptyset\}$ if $\sigma$ is a maximal simplex of $K$ and ${\rm link}_K\sigma=\{\,\}$ if $\sigma\notin K$.
Specifically, if $K=\{\,\}$, then ${\rm  link}_K\sigma=\{\,\}$ for all $\sigma$.
\vspace{3mm}

{\bf Definition 2.1}
For a simplicial complex $K$ on $[m]$
and a sequence of topological (not CW-complex!) pairs $(\underline{X},\underline{A})=\{(X_k,A_k)\}_{k=1}^m$,
the {\it polyhedral product space} ${\cal Z}(K;\underline{X},\underline{A})$
is the subspace of $X_1{\times}{\cdots}{\times}X_m$ defined as follows.
For a subset $\tau$ of $[m]$, define
$$D(\tau)= Y_1{\times}\cdots{\times}Y_m,\quad Y_k=\left\{\begin{array}{cl}
X_k&{\rm if}\,\,k\in \tau, \\
A_k&{\rm if}\,\,k\not\in \tau.
\end{array}
\right.$$
Then ${\cal Z}(K;\underline{X},\underline{A})=\cup_{\tau\in K}\,D(\tau)$.
Empty space $\emptyset$ is allowed in a topological pair and $\emptyset{\times}X=\emptyset$ for all $X$.
Define ${\cal Z}(\{\,\};\underline{X},\underline{A})=\emptyset$.
\vspace{3mm}

Notice that $D(\sigma)={\cal Z}(\Delta\!^\sigma;\underline{X},\underline{A})$,
$D(\emptyset)=A_1{\times}{\cdots}{\times}A_m={\cal Z}(\{\emptyset\};\underline{X},\underline{A})$ and
$D([m])=X_1{\times}{\cdots}{\times}X_m={\cal Z}(\Delta\!^m;\underline{X},\underline{A})$.
But $\emptyset={\cal Z}(\{\,\};\underline{X},\underline{A})$ has no corresponding $D(-)$.
\vspace{3mm}

{\bf Example 2.2} For ${\cal Z}(K;\underline{X},\underline{A})$,
let $S=\{k\,|\,A_k=\emptyset\}$. Then
$${\cal Z}(K;\underline{X},\underline{A})={\cal Z}({\rm link}_KS;\underline{X'},\underline{A'}){\times}(\Pi_{k\in S}X_k),$$
where $(\underline{X'},\underline{A'})=\{(X_k,A_k)\}_{k\notin S}$ and link is as defined in conventions.
\vspace{3mm}

{\bf Definition 2.3} Let $K$ be a simplicial complex with vertex set a subset of $S\neq\emptyset$.
The {\it dual of $K$ relative to} $S$ is the simplicial
complex
$$K^\circ=\{\,S{\setminus}\sigma\,\,|\,\,\sigma\subset S,\,\sigma\!\notin\! K\,\}.$$

It is obvious that
$(K^\circ)^\circ=K$, $(K_1{\cup}K_2)^\circ=(K_1)^\circ{\cap}(K_2)^\circ$ and
$(K_1{\cap}K_2)^\circ=(K_1)^\circ{\cup}(K_2)^\circ$. Specifically, $(\Delta\!^S)^\circ=\{\,\}$
and $(\partial\Delta\!^S)^\circ=\{\emptyset\}$.
\vspace{3mm}

{\bf Theorem 2.4} {\it For ${\cal Z}(K;\underline{X},\underline{A})$, the complement space
$${\cal Z}(K;\underline{X},\underline{A})^c
=(X_1{\times}{\cdots}{\times}X_m){\setminus}{\cal Z}(K;\underline{X},\underline{A})
={\cal Z}(K^\circ;\underline{X},\underline{A}^c),$$
where $(\underline{X},\underline{A}^c)=\{(X_k,A_k^c)\}_{k=1}^m$ with $A_k^c=X_k{\setminus}A_k$
and $K^\circ$ is the dual of $K$ relative to $[m]$. The polyhedral product space ${\cal Z}(K^\circ;\underline{X},\underline{A}^c)$
is called the complement of ${\cal Z}(K;\underline{X},\underline{A})$.
}\vspace{2mm}

{\it Proof}\, For $\sigma\subset[m]$ but $\sigma\neq[m]$ ($\sigma=\emptyset$ is allowed),
$$\begin{array}{l}
\quad(X_1{\times}{\cdots}{\times}X_m)\setminus D(\sigma)\,\,({\rm with\,space\,pair}\,(X_k,A_k))\vspace{2mm}\\
=\cup_{j\notin\sigma}\,X_1{\times}{\cdots}{\times}(X_j{\setminus}A_j){\times}{\cdots}{\times}X_m\vspace{2mm}\\
=\cup_{j\in[m]\setminus\sigma}\,D([m]{\setminus}\{j\})\,\,({\rm with\,space\,pair}\,(X_k,A_k^c))\vspace{2mm}\\
={\cal Z}((\Delta\!^\sigma)^\circ;\underline{X},\underline{A}^c)
\end{array}$$
So for $K\neq\Delta\!^m$ or $\{\,\}$,
$$\begin{array}{l}
\quad{\cal Z}(\Delta\!^m;\underline{X},\underline{A})\setminus{\cal Z}(K;\underline{X},\underline{A})\vspace{2mm}\\
={\cal Z}(\Delta\!^m;\underline{X},\underline{A})\setminus\big(\cup_{\sigma\in K}{\cal Z}(\Delta\!^\sigma;\underline{X},\underline{A})\big)\vspace{2mm}\\
=\cap_{\sigma\in K}\big({\cal Z}(\Delta\!^m;\underline{X},\underline{A})\setminus{\cal Z}(\Delta\!^\sigma;\underline{X},\underline{A})\big)\vspace{2mm}\\
=\cap_{\sigma\in K}{\cal Z}((\Delta\!^\sigma)^\circ;\underline{X},\underline{A}^c)\vspace{2mm}\\
={\cal Z}(K^\circ;\underline{X},\underline{A}^c)\end{array}$$

For $K=\Delta\!^m$ or $\{\,\}$, the above equality holds naturally. \hfill $\Box$
\vspace{3mm}

{\bf Example 2.5} Let $\Bbb F$ be a field and $V$ be a linear space over $\Bbb F$ with base $e_1,{\cdots},e_m$.
For a subset $\sigma=\{i_1,{\cdots},i_s\}\subset[m]$, denote by $\Bbb F(\sigma)$ the subspace of $V$ with base $e_{i_1},{\cdots},e_{i_s}$.
Then for $\Bbb F=\Bbb R$ or $\Bbb C$ and a simplicial complex $K$ on $[m]$, we have
$$V\setminus(\cup_{\sigma\in K}\Bbb R(\sigma))=\Bbb R^m\setminus{\cal Z}(K;\Bbb R,\{0\})={\cal Z}(K^\circ;\Bbb R,\Bbb R{\setminus}\{0\})
\simeq {\cal Z}(K^\circ;D^1,S^0),$$
$$V\setminus(\cup_{\sigma\in K}\Bbb C(\sigma))=\Bbb C^m\setminus{\cal Z}(K;\Bbb C,\{0\})={\cal Z}(K^\circ;\Bbb C,\Bbb C{\setminus}\{0\})
\simeq {\cal Z}(K^\circ;D^2,S^1).$$
This example is applied by Gruji\'{c} and Welkerin in Lemma 2.4 in \cite{GW}.
\vspace{3mm}

{\bf Theorem 2.6} {\it Let $K$ and $K^\circ$ be the dual of each other relative to $[m]$.
The index set $\XX_m=\{(\sigma,\omega)\,|\,\sigma,\omega\subset[m],\,\sigma{\cap}\omega=\emptyset\}$.
For $(\sigma,\omega)\in\XX_m$, define simplicial complex
$K_{\sigma\!,\omega}={\rm link}_K\sigma|_\omega=\{\tau\subset\omega\,|\,\sigma{\cup}\tau\in K\}$
(so $K_{\sigma\!,\omega}=\{\,\}$ if $\sigma\notin K$ or $K=\{\,\}$).
Then for any $(\sigma,\omega)\in\XX_m$ such that $\omega\neq\emptyset$,
$$(K_{\sigma\!,\,\omega})^\circ=(K^\circ)_{\tilde\sigma\!,\,\omega},\quad\w\sigma=[m]{\setminus}(\sigma{\cup}\omega),$$
where $(K_{\sigma\!,\,\omega})^\circ$ is the dual of
$K_{\sigma\!,\,\omega}$ relative to $\omega$.
}\vspace{3mm}

{\it Proof}\, Suppose $\sigma\in K$. Then
$$\begin{array}{l}
\quad (K^\circ)_{\tilde\sigma\!,\,\omega}\vspace{1mm}\\
=\{\eta\,\,|\,\,\eta\subset\omega,\,\,\sigma\!\cup\!(\omega\!\setminus\!\eta)\notin K\}\vspace{1mm}\\
=\{\omega\!\setminus\!\tau\,\,|\,\,\tau\subset\omega,\,\,\sigma\!\cup\!\tau\notin K\}
\quad(\tau=\omega\!\setminus\!\eta)\vspace{1mm}\\
=(K_{\sigma\!,\,\omega})^\circ.
\end{array}$$

If $\sigma\notin K$, then $(K^\circ)_{\tilde\sigma\!,\omega}=\Delta\!^\omega=(K_{\sigma\!,\omega})^\circ$.
\hfill $\Box$\vspace{3mm}

A sequence of simplicial pairs $(\underline{X},\underline{A})=\{(X_k,A_k)\}_{k=1}^m$ in this paper means that
the vertex set of $X_k$ is a subset of $[n_k]$ ($n_k>0$) which is the subset
$$\{s_{k-1}{+}1,s_{k-1}{+}2,\cdots,s_{k-1}{+}n_k\}\quad(s_k=n_1{+}{\cdots}{+}n_k,\,\,s_{0}=0)$$
of $[n]$ with $n=n_1{+}{\cdots}{+}n_m$.

For simplicial complexes $Y_1,{\cdots},Y_m$ such that the vertex set of $Y_k$ is a subset of $[n_k]$,
the union simplicial complex is
$$Y_1*{\cdots}*Y_m=\{\sigma\!\subset\![n]\,|\,\sigma{\cap}[n_k]\in Y_k\,\,{\rm for}\,\,k\!=\!1,{\cdots},m\}.$$

{\bf Definition 2.7} Let $K$ be a simplicial complex on $[m]$ and
$(\underline{X},\underline{A})$ be as above.
The {\it polyhedral product complex} ${\cal S}(K;\underline{X},\underline{A})$
is the simplicial complex on $[n]$ defined as follows. For a subset $\tau\subset[m]$,
define
$$S(\tau)= Y_1*\cdots *Y_m,\quad Y_k=\left\{\begin{array}{cl}
X_k&{\rm if}\,\,k\in \tau, \vspace{1mm}\\
A_k&{\rm if}\,\,k\not\in \tau.
\end{array}
\right.$$
Then ${\cal S}(K;\underline{X},\underline{A})=\cup_{\tau\in K}S(\tau)$.
Void complex $\{\,\}$ is allowed in a simplicial pair and $\{\,\}*X=\{\,\}$ for all $X$.
Define ${\cal Z}(\{\,\};\underline{X},\underline{A})=\{\,\}$.
\vspace{3mm}

{\bf Example 2.8} For ${\cal S}(K;\underline{X},\underline{A})$,
let $S=\{k\,|\,A_k=\{\,\}\}$. Then
$${\cal S}(K;\underline{X},\underline{A})={\cal S}({\rm link}_KS;\underline{X'},\underline{A'}){\times}(\Pi_{k\in S}X_k),$$
where $(\underline{X'},\underline{A'})=\{(X_k,A_k)\}_{k\notin S}$ and link is as defined in conventions.
\vspace{3mm}

{\bf Theorem 2.9} {\it Let ${\cal Z}(K;\underline{Y},\underline{B})$,
$(\underline{Y},\underline{B})=\{(Y_k,B_k)\}_{k=1}^m$ be the polyhedral product space
defined as follows. For each $k$, $(Y_k,B_k)$ is a pair of polyhedral product spaces given by
$$(Y_k,B_k)=\big({\cal Z}(X_k;\underline{U_k},\underline{C_k}),\,{\cal Z}(A_k;\underline{U_k},\underline{C_k})\big),\,\,
(\underline{U_k},\underline{C_k})=\{(U_i,C_i)\}_{i\,=\,s_{k-1}{+}1}^{s_k},$$
where $(X_k,A_k)$ is a simplicial pair on $[n_k]$.
Then
$${\cal Z}(K;\underline{Y},\underline{B})={\cal Z}({\cal S}(K;\underline{X},\underline{A});\underline{U},\underline{C}),$$
where $(\underline{U}{,}\underline{C})=\{(U_k{,}C_k)\}_{k=1}^n$,
$n=n_1{+}{\cdots}{+}n_m$.
\vspace{2mm}

Proof}\, If $K=\{\,\}$ or $U_k=\emptyset$ for some $k$ (this implies $Y_l=\emptyset$ for some $l$), then
${\cal Z}(K;\underline{Y},\underline{B})={\cal Z}({\cal S}(K;\underline{X},\underline{A});\underline{U},\underline{C})=\emptyset$.
So we suppose $K\neq\{\,\}$ and $U_k\neq\emptyset$ for all $k$ in the remaining proof.

We first prove the case $C_k\neq\emptyset$ for all $k$.
If $X_k=\{\,\}$ for some $k$, then $Y_k=\emptyset$ and ${\cal S}(K;\underline{X},\underline{A})=\{\,\}$.
So ${\cal Z}(K;\underline{Y},\underline{B})
={\cal Z}({\cal S}(K;\underline{X},\underline{A});\underline{U},\underline{C})=\emptyset$.
Suppose $X_k\neq\{\,\}$ for all $k$.
Let $S=\{k\,|\,A_k=\{\,\}\}=\{k\,|\,B_k=\emptyset\}$. If $S\not\in K$, then ${\rm link}_KS=\{\,\}$.
From Example~2.2 and Example 2.8 we have
${\cal Z}(K;\underline{Y},\underline{B})={\cal Z}({\cal S}(K;\underline{X},\underline{A});\underline{U},\underline{C})=\emptyset$.
Suppose $S\in K$.
Let $Z^\tau_k=Y_k$ if $k\in\tau$ and $Z^\tau_k=B_k$ if $k\notin\tau$.
For $\tau_k\subset[n_k]$, $W^{\tau_k}_t=U_k$ if $t\in\tau_k$ and $W^{\tau_k}_t=C_k$ if $k\in[n_k]{\setminus}\tau_k$. Then
$$\begin{array}{l}
\quad {\cal Z}(K;\underline{Y},\underline{B})\vspace{2mm}\\
=\cup_{\tau\in K} Z^\tau_1\times{\cdots}\times Z^\tau_m\vspace{2mm}\\
=\cup_{\tau,\tau_1,{\cdots},\tau_m}\,(W^{\tau_1}_1{\times}{\cdots}{\times}W^{\tau_1}_{n_1}){\times}{\cdots}{\times}
(W^{\tau_m}_{s_{m-1}+1}{\times}{\cdots}{\times}W^{\tau_m}_n)\vspace{2mm}\\
=\cup_{(\tau_1,{\cdots},\tau_m)\in{\cal S}(K;\underline{X},\underline{A})}\,(W^{\tau_1}_1{\times}{\cdots}{\times}W^{\tau_1}_{n_1}){\times}{\cdots}{\times}
(W^{\tau_m}_{s_{m-1}+1}{\times}{\cdots}{\times}W^{\tau_m}_n)\vspace{2mm}\\
={\cal Z}({\cal S}(K;\underline{X},\underline{A});\underline{U},\underline{C}),
\end{array}$$
where $\tau,\tau_1,{\cdots},\tau_m$ are taken over all subsets such that $\tau\in K$ and $S\subset\tau$, $\tau_k\in X_k$ if $k\in\tau$
and $\tau_k\in A_k$ if $k\notin\tau$.

Now we prove the case $\sigma=\{k\,|\, C_k=\emptyset\}\neq\emptyset$.
Let $\sigma_k=\sigma{\cap}[n_k]$. Then from Example 2.2 we have
$$\begin{array}{l}
\,Y_k=Y'_k\times(\times_{k\in\sigma_k}U_k),\quad
\,Y'_k={\cal Z}({\rm link}_{X_k}\sigma_k;\underline{U'_k},\underline{C'_k}),\vspace{2mm}\\
B_k=B'_k\times(\times_{k\in\sigma_k}U_k),\quad
B'_k={\cal Z}({\rm link}_{A_k}\sigma_k;\underline{U'_k},\underline{C'_k}),\vspace{2mm}\\
{\cal Z}(K;\underline{Y},\underline{B})={\cal Z}(K;\underline{Y'},\underline{B'})\times(\times_{k\in\sigma}U_k),
\end{array}$$
where $(\underline{U'_k},\underline{C'_k})=\{(U_k,C_k)\}_{k\in[n_k]\setminus\sigma_k}$
and $(\underline{Y'},\underline{B'})=\{(Y'_k,B'_k)\}_{k\notin\sigma}$.
Denote  ${\rm link}_{(\underline{X},\underline{A})}\sigma=\{({\rm link}_{X_k}\sigma_k,{\rm link}_{A_k}\sigma_k)\}_{k=1}^m$.
Then ${\cal S}(K;{\rm link}_{(\underline{X},\underline{A})}\sigma)={\rm link}_{{\cal S}(K;\underline{X},\underline{A})}\sigma$
(this equality is a special case of Theorem 2.10 for $(\sigma\!,\,\omega)=(\sigma,[n]{\setminus}\sigma)$ and so
$(-)_{\sigma\!,\,\omega}={\rm link}_{(-)}\sigma$). So
$$\begin{array}{l}
\quad{\cal Z}(K;\underline{Y},\underline{B})\vspace{2mm}\\
={\cal Z}(K;\underline{Y'},\underline{B'})\times(\times_{k\in\sigma}U_k)\vspace{2mm}\\
={\cal Z}({\cal S}(K;{\rm link}_{(\underline{X},\underline{A})}\sigma);\underline{U'},\underline{C'})\times(\times_{k\in\sigma}U_k)\vspace{2mm}\\
={\cal Z}({\rm link}_{{\cal S}(K;(\underline{X},\underline{A})}\sigma;\underline{U'},\underline{C'})\times(\times_{k\in\sigma}U_k)\vspace{2mm}\\
={\cal Z}({\cal S}(K;
\underline{X},\underline{A});\underline{U},\underline{C}),
\end{array}$$
where $(\underline{U'},\underline{C'})=\{(U_k,C_k)\}_{k\notin\sigma}$.
\hfill$\Box$\vspace{3mm}

With this theorem we see that to compute the cohomology algebra of ${\cal Z}(K;\underline{Y},\underline{B})$,
we have to compute the universal algebra of ${\cal S}(K;\underline{X},\underline{A})$, which is the central work
of this paper.
\vspace{3mm}

{\bf Theorem 2.10} {\it \, For the ${\cal S}(K;\underline{X},\underline{A})$ in Definition 2.7 and
$(\sigma,\omega)\in\XX_n$ ( the simplicial complex $(-)_{\sigma,\omega}$ is as defined in Theorem 2.6),
$${\cal S}(K;\underline{X},\underline{A})_{\sigma,\,\omega}
={\cal S}(K;\underline{X}\,_{\sigma,\,\omega},\underline{A}\,_{\sigma,\,\omega}),$$
where $(\underline{X}\,_{\sigma,\,\omega},\underline{A}\,_{\sigma,\,\omega})\!=\!
\{(\,(X_k)_{\sigma_k,\,\omega_k},(A_k)_{\sigma_k,\,\omega_k})\}_{k=1}^m$,
$\sigma_k\!=\!\sigma{\cap}[n_k]$, $\omega_k\!=\!\omega{\cap}[n_k]$.
\vspace{2mm}

Proof}\, By definition, $(K{\cup}L)_{\sigma\!,\,\omega}=K_{\sigma\!,\,\omega}{\cup}L_{\sigma\!,\,\omega}$
and $(Y_1*{\cdots}*Y_m)_{\sigma\!,\,\omega}=(Y_1)_{\sigma_1,\,\omega_1}*{\cdots}*(Y_m)_{\sigma_m,\,\omega_m}$.
Let $Y^\tau_k=X_k$ if $k\in\tau$ and $Y^\tau_k=A_k$ if $k\notin\tau$ (the void complex $\{\,\}$ is allowed).
Then
\begin{eqnarray*}&&{\cal S}(K;\underline{X},\underline{A})_{\sigma\!,\,\omega}\vspace{1mm}\\
&=&\cup_{\tau\in K}(Y^\tau_1*{\cdots}*Y^\tau_m)_{\sigma\!,\,\omega}\vspace{1mm}\\
&=&\cup_{\tau\in K}(Y^\tau_1)_{\sigma_1,\,\omega_1}*\cdots*(Y^\tau_m)_{\sigma_m,\,\omega_m}\vspace{1mm}\\
&=&{\cal S}(K;\underline{X}\,_{\sigma,\,\omega},\underline{A}\,_{\sigma,\,\omega}).
\end{eqnarray*}
\hfill$\Box$

Notice that the dual of ${\cal S}(K;\underline{X},\underline{A})$ relative to $[n]$ is in general not a polyhedral product complex.
But  the following special type of polyhedral product complexes is closed with respect to duality.
\vspace{3mm}

{\bf Definition 2.11} ${\cal S}(K;L_1,{\cdots},L_m)$ is the polyhedral product complex ${\cal S}(K;\underline{X},\underline{A})$
such that each pair $(X_k,A_k)=(\Delta\!^{n_k},L_k)$.
\vspace{3mm}

The complex ${\cal S}(K;K_1,{\cdots},K_m)$ is the composition complex $K(K_1,{\cdots},K_m)$ in Definition 4.5 of \cite{AY}.
\vspace{3mm}

{\bf Theorem 2.12} {\it Let ${\cal S}(K;L_1,{\cdots},L_m)^\circ$ be the dual of ${\cal S}(K;L_1,{\cdots},L_m)$ relative to $[n]$.
Then
$${\cal S}(K;L_1,{\cdots},L_m)^\circ={\cal S}(K^\circ;L_1^\circ,{\cdots},L_m^\circ),$$
where $K^\circ$ is the dual of $K$ relative to $[m]$ and $L_k^\circ$ is the dual  of $L_k$ relative to $[n_k]$.
So if $K$ and all $L_k$ are self dual ($X=X^\circ$ relative to its non-empty vertex set),
then ${\cal S}(K;L_1,{\cdots},L_m)$ is self dual.
\vspace{2mm}

Proof}\, For $\sigma\subset[m]$ but $\sigma\neq[m]$ ($\sigma=\emptyset$, $L_k=\Delta\!^{n_k}$ or $\{\,\}$ are allowed),
$$\begin{array}{l}
\quad {\cal S}(\Delta\!^\sigma;L_1,{\cdots},L_m)^\circ\vspace{2mm}\\
=\{[n]\!\setminus\!\tau\,|\,\tau\in\cup_{j\notin\sigma}\Delta\!^{n_1}*{\cdots}*(\Delta\!^{n_j}{\setminus}L_j)*{\cdots}*\Delta\!^{n_m}\}\vspace{2mm}\\
=\cup_{j\notin\sigma}\,\Delta\!^{n_1}*{\cdots}*L_j^\circ*{\cdots}*\Delta\!^{n_m}\vspace{2mm}\\
={\cal S}((\Delta\!^\sigma)^\circ;L_1^\circ,{\cdots},L_m^\circ),
\end{array}$$
So for $K\neq[m]$ or $\{\,\}$,
$$\begin{array}{l}
\quad {\cal S}(K;L_1,{\cdots},L_m)^\circ\vspace{2mm}\\
=(\cup_{\sigma\in K}{\cal S}((\Delta\!^\sigma);L_1,{\cdots},L_m))^\circ\vspace{2mm}\\
={\cal S}(\cap_{\sigma\in K}(\Delta\!^\sigma)^\circ;L_1^\circ,{\cdots},L_m^\circ)\vspace{2mm}\\
={\cal S}(K^\circ;L_1^\circ,{\cdots},L_m^\circ).
\end{array}$$

For $K=\Delta\!^m$ or $\{\,\}$, the equality holds naturally.
\hfill$\Box$\vspace{3mm}

\section{Homology and Cohomology Group}

This is a paper following \cite{Z}.
All the basic definitions such as indexed groups and (co)chain complexes, diagonal tensor product,
etc., are as in \cite{Z}.

In this section, we compute the reduced simplicial (co)homology group and
the (right) total (co)homology group of polyhedral product complexes uniformly from the point of view of split inclusion.
\vspace{3mm}

{\bf Conventions} In this paper, a group $A_*^\Lambda=\oplus_{\alpha\in\Lambda}\,A_*^\alpha$ indexed by $\Lambda$
is simply denoted by $A_*$ when there is no confusion.
So is the (co)chain complex case.
The diagonal tensor product $A_*^\Lambda{\otimes}_\Lambda B_*^\Lambda$ in \cite{Z} is simply denoted by
$A_*^\Lambda\,\widehat\otimes\,B_*^\Lambda$ in this paper (so $\Lambda$ can not be abbreviated in this case).
\vspace{3mm}

{\bf Definition 3.1} Let $A_*=\oplus_{\alpha\in\Lambda}A_*^\alpha$, $B_*=\oplus_{\alpha\in\Lambda}B_*^\alpha$
be two groups indexed by the same set $\Lambda$.
An {\it indexed group homomorphism} $f\colon A_*\to B_*$ is the direct sum
$f=\oplus_{\alpha\in\Lambda}f_\alpha$ such that each $f_\alpha\colon A_*^\alpha\to B_*^{\alpha}$
is a graded group homomorphism. Define groups indexed by $\Lambda$ as follows.
$$\begin{array}{ll}
{\rm ker}\,f=\oplus_{\alpha\in\Lambda}{\rm ker}\,f_\alpha,&{\rm coker}\,f=\oplus_{\alpha\in\Lambda}{\rm coker}\,f_\alpha,\vspace{1mm}\\
{\rm im}\,f\,=\oplus_{\alpha\in\Lambda}\,{\rm im}\,f_\alpha,&{\rm coim}\,f\,=\oplus_{\alpha\in\Lambda}\,{\rm coim}\,f_\alpha.
\end{array}$$

For indexed group homomorphism $f=\oplus_{\alpha\in\Lambda}\,f_\alpha$ and $g=\oplus_{\beta\in\Gamma}\,g_\beta$,
their tensor product $f{\otimes}g$ is naturally an indexed group homomorphism with
$f{\otimes}g=\oplus_{(\alpha,\beta)\in\Lambda{\times}\Gamma}f_\alpha{\otimes}g_\beta$.

For indexed group homomorphism $f=\oplus_{\alpha\in\Lambda}\,f_\alpha$ and $g=\oplus_{\alpha\in\Lambda}\,g_\alpha$ indexed by
the same set,
their diagonal tensor product $f\widehat\otimes g$ is the indexed group homomorphism
$f\widehat\otimes g=\oplus_{\alpha\in\Lambda}f_\alpha{\otimes}g_\alpha$.

Similarly, we have the definition of {\it indexed (co)chain homomorphism}
by replacing the indexed groups in the above definition by indexed (co)chain complexes.
\vspace{3mm}

{\bf Definition 3.2} An indexed group homomorphism $\theta\colon U_*\to V_*$
is called {\it a split homomorphism} if ${\rm ker}\,\theta$, ${\rm coker}\,\theta$ and ${\rm im}\,\theta$ are all free groups.

An indexed chain homomorphism
$\vartheta\colon (C_*,d)\to(D_*,d)$ with induced homology group homomorphism $\theta\colon U_*\to V_*$
($U_*\!=\!H_*(C_*)$, $V_*\!=\!H_*(D_*)$) is called a {\it split inclusion} if $C_*$ is a chain subcomplex of the free complex $D_*$ and
$\theta$ is a split homomorphism.
\vspace{3mm}

For a topological pair $(X,A)$, let $\vartheta\colon S_*(A)\to S_*(X)$ be the singular chain complex inclusion.
Regard this inclusion as an indexed chain homomorphism such that the index set has only one element.
Then $\vartheta$ is a split inclusion if and only if $(X,A)$ is homology split as defined in Definition 2.1 in \cite{Z}.
The work of this and the next section is just to generalize all the work in [17] from
the singular chain complex case to indexed total chain complex case.
\vspace{3mm}

{\bf Definition 3.3} Let $\theta\colon U_*\to V_*$ be a split homomorphism
with dual homomorphism $\theta^\circ\colon V^*\to U^*$.
The index set $\XX=\XX_1$ ($\XX_m$ is as defined in Definition 2.6) and $\RR=\{\,\{\emptyset,\emptyset\},\{\emptyset,\{1\}\}\,\}\subset\XX$.
$\SS=\XX$ or $\RR$.

The indexed groups $H_*^{\SS}(\theta)=\oplus_{s\in\SS}\,H_*^{s}(\theta)$
and its dual groups $H^{\,*}_{\SS}(\theta^\circ)=\oplus_{s\in\SS}\,H^*_{s}(\theta^\circ)$
are given by
$$H_*^s(\theta)=\left\{\begin{array}{ll}
{\rm coker}\,\theta&{\rm if}\,\,s={\scriptstyle\{\{1\},\emptyset\}},\\
{\rm ker}\,\theta&{\rm if}\,\,s={\scriptstyle\{\emptyset,\{1\}\}},\\
{\rm im}\,\theta&{\rm if}\,\,s={\scriptstyle\{\emptyset,\emptyset\}},
\end{array}\right.$$
$$H^*_s(\theta^\circ)=\left\{\begin{array}{ll}
{\rm ker}\,\theta^\circ&{\rm if}\,\,s={\scriptstyle\{\{1\},\emptyset\}},\\
{\rm coker}\,\theta^\circ&{\rm if}\,\,s={\scriptstyle\{\emptyset,\{1\}\}},\\
{\rm im}\,\theta^\circ&{\rm if}\,\,s={\scriptstyle\{\emptyset,\emptyset\}},
\end{array}\right.
$$

The indexed chain complexes $(C_*^{\SS}(\theta),d)=\oplus_{s\in\SS}\,(C_*^{s}(\theta),d)$
and its dual cochain complexes $(C^{\,*}_{\!\SS}(\theta^\circ),\delta)=\oplus_{s\in\SS}\,(C^*_{s}(\theta^\circ),\delta)$
are given by
$$C_*^s(\theta)=\left\{\begin{array}{ll}
{\rm coker}\,\theta&{\rm if}\,\,s={\scriptstyle\{\{1\},\emptyset\}},\\
{\rm ker}\,\theta{\oplus}\Sigma{\rm ker}\,\theta&{\rm if}\,\,s={\scriptstyle\{\emptyset,\{1\}\}},\\
{\rm im}\,\theta&{\rm if}\,\,s={\scriptstyle\{\emptyset,\emptyset\}},
\end{array}\right.\quad\quad$$
$$C^*_s(\theta^\circ)=\left\{\begin{array}{ll}
{\rm ker}\,\theta^\circ&{\rm if}\,\,s={\scriptstyle\{\{1\},\emptyset\}},\\
{\rm coker}\,\theta^\circ{\oplus}\Sigma{\rm coker}\,\theta^\circ&{\rm if}\,\,s={\scriptstyle\{\emptyset,\{1\}\}},\\
{\rm im}\,\theta^\circ&{\rm if}\,\,s={\scriptstyle\{\emptyset,\emptyset\}},
\end{array}\right.$$
where $d$ is trivial on $C_*^{\emptyset,\emptyset}(\theta)$ and $C_*^{\{1\},\emptyset}(\theta)$
and is the desuspension isomorphism on $C_*^{\emptyset,\{1\}}(\theta)$.

Let $\theta=\oplus_{\alpha\in\Lambda}\theta_{\alpha}$.
Then $H_*^{\SS}\!(\theta)$ is also a group indexed by $\Lambda$
and so denoted by $H_*^{\SS}\!(\theta)=H_*^{\SS\!;\Lambda}(\theta)=\oplus_{s\in\SS\!,\alpha\in\Lambda}\,
H_*^{s;\alpha}(\theta)$ with
$$H_*^{s;\alpha}(\theta)=\left\{\begin{array}{ll}
{\rm coker}\,\theta_\alpha&{\rm if}\,\,s={\scriptstyle\{\{1\},\emptyset\}},\\
{\rm ker}\,\theta_\alpha&{\rm if}\,\,s={\scriptstyle\{\emptyset,\{1\}\}},\\
{\rm im}\,\theta_\alpha&{\rm if}\,\,s={\scriptstyle\{\emptyset,\emptyset\}}.
\end{array}\right.$$
Other cases are similar.
\vspace{3mm}

{\bf Theorem 3.4}\, {\it For a split inclusion $\vartheta\colon(C_*,d)\to(D_*,d)$
with induced homology homomorphism $\theta\colon U_*\to V_*$,
there are quotient chain homotopy equivalences  $q$ and $q'$
satisfying the following commutative diagram
$$\begin{array}{ccc}
(C_*,d)&\stackrel{q'}{\longrightarrow}&(U_*,d)\,\vspace{1mm}\\
^\vartheta\downarrow\,\,&&^{\vartheta'}\downarrow\quad\\
(D_*,d)&\stackrel{q}{\longrightarrow}&(C_*^{\XX}(\theta),d),
\end{array}
$$
where $\vartheta'$ is the inclusion by identifying $U_*={\rm ker}\,\theta\oplus{\rm coim}\,\theta$ with
${\rm ker}\,\theta\oplus{\rm im}\,\theta\subset C_*^{\XX}(\theta)$ ($d$ is trivial on $U_*$).

There are also isomorphisms $\phi$ and $\phi'$ of chain complexes indexed by $\XX$
satisfying the following commutative diagram
$$\begin{array}{ccc}
(U_*,d)&\stackrel{\phi'}{\longrightarrow}&S_*^\XX\,\widehat\otimes\, H_*^{\XX}\!(\theta)\vspace{1mm}\\
^{\vartheta'}\downarrow\,\,&&^{i\,\widehat\otimes\, 1}\downarrow\quad\quad\quad\\
(C_*^{\XX}\!(\theta),d)&\stackrel{\phi}{\longrightarrow}&(T_*^\XX\,\widehat\otimes\, H_*^{\XX}\!(\theta),d)
\end{array}
$$
where $S_*^\XX$ and $T_*^\XX$ are as in Definition 4.2 and Theorem 4.3 in} \cite{Z}, {\it $1$ is the identity and $i$ is the inclusion.

If $\theta$ is an epimorphism, then
$H_*^{\XX}\!(\theta)=H_*^{\RR}(\theta)=U_*$ by identify ${\rm im}\,\theta$
with ${\rm coim}\,\theta$ and so  all $\XX$ is replaced by $\RR$.
\vspace{2mm}

Proof}\, Take a representative $a_i$ in $C_*$ for every generator of ${\rm ker}\,\theta$
and let $\overline a _i\in D_*$ be any element such that $d\overline a _i= a_i$.
Take a representative $b_j$ in $C_*$ for every generator of ${\rm im}\,\theta$.
Take a representative $c_k$ in $D_*$ for every generator of ${\rm coker}\,\theta$.
So we may regard $U_*$ as the chain subcomplex of $C_*$ freely generated by all $a_i$'s and $b_j$'s
and regard $(C_*^{\XX\!}(\theta),d)$ as the chain subcomplex of $D_*$ freely generated by all
$a_i$'s, $\overline a _i$'s, $b_j$'s and $c_k$'s.
Then we have the following commutative diagram of short exact sequences of chain complexes
$$\begin{array}{ccccccc}
 0\to& U_*&\stackrel{i}{\longrightarrow}& C_*&\stackrel{j}{\longrightarrow}&C_*/U_*&\to0\vspace{1mm}\\
 &^{\vartheta'}\downarrow\quad&&^\vartheta\downarrow\,\quad&&\downarrow&\\
0\to& C_*^{\XX\!}(\theta)&\stackrel{i}{\longrightarrow}& D_*&\stackrel{j}{\longrightarrow}&D_*/C_*^{\XX\!}(\theta)&\to0.
  \end{array}
$$
Since all the complexes are free, $i$'s have group homomorphism inverse. $H_*(C_*/U_*)=0$ and $H_*(D_*/C_*^{\XX\!}(\theta))=0$
imply that the inverse of $i$'s are complex homomorphisms. So we may take $q,q'$ to be the inverse of $i$'s.

$\phi$ is defined as shown in the following table.
\begin{center}
\begin{tabular}{|c|c|c|c|c|}
\hline
{\rule[-2mm]{0mm}{6mm}$\quad x\in$}&${\rm coker}\,\theta$&$\Sigma\,{\rm ker}\,\theta$
&${\rm ker}\,\theta$&${\rm im}\,\theta$\\
\hline
{\rule[-2mm]{0mm}{7mm}$\phi(x)=$}&$\alpha\,\widehat\otimes\,x$&$\beta\,\widehat\otimes\,dx$
&$\gamma\,\widehat\otimes\,x$&$\eta\,\widehat\otimes\,x$\\
\hline
\end{tabular}
\end{center}
\hfill$\Box$\vspace{3mm}

{\bf Definition 3.5} For $k=1,{\cdots},m$,
let $\vartheta_k\colon\! ((C_k)_*,d)\to((D_k)_*,d)$ be a split inclusion with induced
homology group homomorphism $\theta_k\colon\! (U_k)_*\to (V_k)_*$.
Denote $\underline{\vartheta}=\{\vartheta_k\}_{k=1}^m$, $\underline{\theta}=\{\theta_k\}_{k=1}^m$
and their dual $\underline{\vartheta}^\circ=\{\vartheta_k^\circ\}_{k=1}^m$, $\underline{\theta}^\circ=\{\theta_k^\circ\}_{k=1}^m$.
The index set $\SS=\XX$ or $\RR$.

The indexed group $H_*^{\SS_m}(\underline{\theta})$ and its dual group $H^{\,*}_{\!\SS_m}(\underline{\theta}^\circ)$
(the index set $\SS_m$ as in Theorem 2.6) are given by
$$H_*^{\SS_m}(\underline{\theta})=
H_*^{\SS}(\theta_1)\otimes{\cdots}\otimes H_*^{\SS}(\theta_m),\,\,
H^{\,*}_{\!\SS_m}(\underline{\theta}^\circ)=
H^{\,*}_{\!\SS}(\theta_1^\circ)\otimes{\cdots}\otimes H^{\,*}_{\!\SS}(\theta_m^\circ).$$
Denote $H_*^{\SS_m}(\underline{\theta})=\oplus_{(\sigma\!,\,\omega)\in\SS_m}H_*^{\sigma\!,\,\omega}(\underline{\theta})$,
$H^{\,*}_{\!\SS_m}(\underline{\theta}^\circ)=\oplus_{(\sigma\!,\,\omega)\in\SS_m}H^*_{\sigma\!,\,\omega}(\underline{\theta}^\circ)$. Then
$$H_*^{\sigma\!,\,\omega}(\underline{\theta})
= H_1{\otimes}{\cdots}{\otimes}H_m,\quad
H_k=\left\{\begin{array}{cl}
H_*^{\{1\},\emptyset}(\theta_{k})&{\rm if}\, k\in\sigma,\vspace{1mm}\\
H_*^{\emptyset,\{1\}}(\theta_{k})&{\rm if}\, k\in\omega,\vspace{1mm}\\
H_*^{\emptyset,\emptyset}(\theta_{k})&{\rm otherwise},
\end{array}\right.$$
$$H^*_{\sigma\!,\,\omega}(\underline{\theta}^\circ)
= H^1{\otimes}{\cdots}{\otimes}H^m,\quad
H^k=\left\{\begin{array}{cl}
H_*^{\{1\},\emptyset}(\theta_{k}^\circ)&{\rm if}\, k\in\sigma,\vspace{1mm}\\
H_*^{\emptyset,\{1\}}(\theta_{k}^\circ)&{\rm if}\, k\in\omega,\vspace{1mm}\\
H_*^{\emptyset,\emptyset}(\theta_{k}^\circ)&{\rm otherwise}.
\end{array}\right.$$

The indexed chain complex $(C_*^{\SS_m}(\underline{\theta}),d)$ and its dual cochain complex
$(C^{\,*}_{\!\SS_m}(\underline{\theta}^\circ),\delta)$ are given by
$$C_*^{\SS_m}(\underline{\theta})=C_*^{\SS}\!(\theta_1)\otimes{\cdots}\otimes C_*^{\SS}\!(\theta_m),\,\,
C^{\,*}_{\!\SS_m}(\underline{\theta}^\circ)=C^{\,*}_{\!\SS}(\theta_1^\circ)\otimes{\cdots}\otimes C^{\,*}_{\!\SS}(\theta_m^\circ).
\vspace{3mm}$$

{\bf Definition 3.6} Let $K$ be a simplicial complex on $[m]$
and everything else be as in Definition 3.5.

The indexed chain complex $(C_*^{\SS_m}(K;\underline{\theta}),d)$ is the subcomplex of
$(C_*^{\SS_m}(\underline{\theta}),d)$ defined as follows.
For a subset $\tau$ of $[m]$, define
$$(H_*(\tau),d)=(H_1{\otimes}\cdots{\otimes}H_m,d),\quad (H_k,d)=\left\{\begin{array}{cl}
(C_*^{\SS}\!(\theta_k),d)&{\rm if}\,\,k\in \tau,\vspace{1mm}\\
((U_k)_*,d)&{\rm if}\,\,k\not\in \tau.
\end{array}\right.$$
Then $(C_*^{\SS_m}(K;\underline{\theta}),d)=(+_{\tau\in K}\,H_*(\tau),d)$. Define
$(C_*^{\SS_m}(\{\,\};\underline{\theta}),d)=0$.

So the dual cochain complex $(C^{\,*}_{\!\SS_m}(K;\underline{\theta}^\circ),\delta)$ of
$(C_*^{\SS_m}(K;\underline{\theta}),d)$ is a quotient complex of $(C^{\,*}_{\!\SS_m}(\underline{\theta}^\circ),\delta)$.

The chain complex $(C_*(K;\underline{\vartheta}),d)$ is the subcomplex of $((D_1)_*{\otimes}{\cdots}{\otimes}(D_m)_*,d)$
defined as follows.
For a subset $\tau$ of $[m]$, define
$$(E_*(\tau),d)=(E_1{\otimes}\cdots{\otimes}E_m,d),\quad (E_k,d)=\left\{\begin{array}{cl}
((D_k)_*,d)&{\rm if}\,\,k\in \tau,\vspace{1mm}\\
((C_k)_*,d)&{\rm if}\,\,k\not\in \tau.
\end{array}
\right.$$
Then $(C_*(K;\underline{\vartheta}),d)=(+_{\tau\in K}\,E_*(\tau),d)$.
Define $(C_*(\{\,\};\underline{\vartheta}),d)=0$.

So the  dual cochain complex $(C^*(K;\underline{\vartheta}^\circ),\delta)$ of $(C_*(K;\underline{\vartheta}),d)$
is a quotient complex of $((D_1)^*{\otimes}{\cdots}{\otimes}(D_m)^*,\delta)$.
\vspace{3mm}

{\bf Theorem 3.7}\, {\it For the $K$, $\underline{\vartheta}$ and $\underline{\theta}$ in Definition 3.5 and Definition 3.6,
there is a quotient chain homotopy equivalence (\,$\XX_m$ neglected)
$$\varphi_{(K;\underline{\vartheta})}\colon(C_*(K;\underline{\vartheta}),d)
\stackrel{\simeq}{\longrightarrow}
(C_*^{\XX_m}(K;\underline{\theta}),d)$$
and an isomorphism of chain complexes indexed by $\XX_m$
$$\phi_{(K;\underline{\theta})}\colon(C_*^{\XX_m}(K;\underline{\theta}),d)
\stackrel{\cong}{\longrightarrow}
(T_*^{\XX_m}(K)\,\widehat\otimes\,H_*^{\XX_m}(\underline{\theta}),d).$$
So we have (co)homology group isomorphisms
$$H_*(C_*(K;\underline{\vartheta}))\cong
H_*^{\XX_m}(K)\,\widehat\otimes\,H_*^{\XX_m}(\underline{\theta})
=H_*^{\XX_m}(K)\,\widehat\otimes\,\big(H_*^{\XX}\!(\theta_1){\otimes}{\cdots}{\otimes}H_*^{\XX}\!(\theta_m)\big),$$
$$H^*(C^*(K;\underline{\vartheta}^\circ))\cong
H^{\,*}_{\!\XX_m}(K)\,\widehat\otimes\,H^{\,*}_{\!\XX_m}(\underline{\theta}^\circ)
=H^*_{\!\XX_m}(K)\,\widehat\otimes\,\big(H^*_{\!\XX}(\theta_1^\circ){\otimes}{\cdots}{\otimes}H^*_{\!\XX}(\theta_m^\circ)\big),$$
where $T_*^{\XX_m}(K)$, $T^{\,*}_{\!\XX_m}(K)$, $H_*^{\XX_m}(K)$ and $H^{\,*}_{\!\XX_m}(K)$ are as in Definition 4.5 and Theorem 4.7 in {\rm \cite{Z}}.

If each $\theta_k$ is an epimorphism, then
$H_*^{\XX}\!(\theta_k)=H_*^{\RR}(\theta_k)=(U_k)_*$, $H^{\,*}_{\XX}(\theta_k^\circ)=H^{\,*}_{\RR}(\theta_k^\circ)=(U_k)^*$
and so all $\XX$ is replaced by $\RR$.
\vspace{2mm}

Proof}\, Denote by $q_k$ and $q'_k$ the chain homotopy equivalence $q$ and $q'$ in Theorem 3.4 for
$\vartheta=\vartheta_k$.
For $\sigma\subset[m]$, let $\varphi_\sigma=p_1{\otimes}{\cdots}{\otimes}p_m$,
where $p_k=q_k$ if $k\in\sigma$ and $p_k=q'_k$ if $k\notin\sigma$.
So $\varphi_\sigma$ is a chain homotopy equivalence.
Then $\varphi_{(K;\underline{\vartheta})}=+_{\sigma\in K}\,\varphi_\sigma$ is also a chain homotopy equivalence.

Denote by $\phi_k$ and $\phi'_k$ the isomorphism $\phi$ and $\phi'$ in Theorem 3.4 for
$\theta=\theta_k$.
For $\sigma\subset[m]$, let $\phi_\sigma=r_1{\otimes}{\cdots}{\otimes}r_m$,
where $r_k=\phi_k$ if $k\in\sigma$ and $r_k=\phi'_k$ if $k\notin\sigma$.
So $\phi_\sigma$ is an isomorphism.
Then $\phi_{(K;\underline{\theta})}=+_{\sigma\in K}\,\phi_\sigma$ is also an isomorphism.
\hfill$\Box$\vspace{3mm}

{\bf Definition 3.8} A polyhedral product complex ${\cal S}(K;\underline{X},\underline{A})$ is {\it homology split}
if the reduced simplicial homology homomorphism $$\iota_k\colon\w H_*(A_k)\to\w H_*(X_k)$$ induced by inclusion is split
for $k=1,{\cdots},m$.

A polyhedral product complex ${\cal S}(K;\underline{X},\underline{A})$ is {\it total homology split}
if the reduced simplicial homology homomorphism
$$\iota_{\sigma_k,\omega_k}\colon \w H_{*}((A_k)_{\sigma_k,\omega_k})
\to \w H_{*}((X_k)_{\sigma_k,\omega_k})\quad((-)_{\sigma\!,\,\omega}\,\,{\rm as\,\,in\,\,Theorem\, 2.6})$$ induced by inclusion is split
for $k=1,{\cdots},m$ and all $(\sigma_k,\omega_k)\in\XX_{n_k}$.
\vspace{3mm}

{\bf Theorem 3.9}{\it\, For homology split ${\cal S}(K;\underline{X},\underline{A})$,
$$\begin{array}{l}
\quad\w H_{*-1}({\cal S}(K;\underline{X},\underline{A}))\vspace{2mm}\\
\cong H_*^{\XX_m}(K)\,\widehat\otimes\,H_*^{\XX_m}(\underline{X},\underline{A})\vspace{2mm}\\
=H_*^{\XX_m}(K)\,\widehat\otimes\,\big(H_*^{\XX}\!(X_1,A_1){\otimes}{\cdots}{\otimes}H_*^{\XX}\!(X_m,A_m)\big),
\end{array}$$
$$\begin{array}{l}
\quad\w H^{*-1}({\cal S}(K;\underline{X},\underline{A}))\vspace{2mm}\\
\cong H^*_{\XX_m}(K)\,\widehat\otimes\,H^*_{\XX_m}(\underline{X},\underline{A})\vspace{2mm}\\
=H_*^{\XX_m}(K)\,\widehat\otimes\,\big(H^{\,*}_{\!\XX}(X_1,A_1){\otimes}{\cdots}{\otimes}H^{\,*}_{\!\XX}(X_m,A_m)\big),
\end{array}$$
where $H_*^{\XX_m}(-)=\oplus_{(\sigma,\omega)\in\XX_m}\,H_*^{\sigma\!,\,\omega}(-)$,
$H^*_{\XX_m}(-)=\oplus_{(\sigma,\omega)\in\XX_m}\,H^*_{\sigma\!,\,\omega}(-)$ with
$$H_*^{\sigma\!,\omega}(K)=\w H_{*-1}(K_{\sigma\!,\omega}),\quad H^*_{\sigma\!,\omega}(K)=\w H^{*-1}(K_{\sigma\!,\omega}),$$
$$H_*^{\XX}\!(X_k,A_k)=H_{*-1}^\XX(\iota_k),\quad H^*_{\!\XX}(X_k,A_k)=H^{*-1}_{\!\XX}(\iota_k^\circ),$$
$$H_*^{\sigma\!,\,\omega}(\underline{X},\underline{A})
=H_1{\otimes}{\cdots}{\otimes}H_m,\quad
H_k=\left\{\begin{array}{rl}
\Sigma{\rm coker}\,\iota_k&{\rm if}\, k\in\sigma,\vspace{1mm}\\
\Sigma{\rm ker}\,\iota_k&{\rm if}\, k\in\omega,\vspace{1mm}\\
\Sigma{\rm im}\,\iota_k&{\rm otherwise},
\end{array}\right.$$
$$H^*_{\sigma\!,\,\omega}(\underline{X},\underline{A})
=H^1{\otimes}{\cdots}{\otimes}H^m,\quad
H^k=\left\{\begin{array}{rl}
\Sigma{\rm ker}\,\iota_k^\circ&{\rm if}\, k\in\sigma,\vspace{1mm}\\
\Sigma{\rm coker}\,\iota_k^\circ&{\rm if}\, k\in\omega,\vspace{1mm}\\
\Sigma{\rm im}\,\iota_k^\circ&{\rm otherwise},
\end{array}\right.$$
where $\iota_k$ is as in Definition 3.8 with dual $\iota_k^\circ$ and
$\Sigma$ means suspension.

If each $\iota_k$ is an epimorphism, then all $\XX$ is replaced by $\RR$ and we have $H_*^\RR(X_k,A_k)=H_*(A_k)$.

If the reduced simplicial (co)homology is taken over a field, then the conclusion holds for all
polyhedral product complexes.
\vspace{2mm}

Proof}\, A corollary of Theorem 3.7 by
taking $\underline{\vartheta}=\{\vartheta_k\}_{k=1}^m$ with split inclusion
$\vartheta_k\colon(\Sigma\w C_*(A_k),d)\to(\Sigma\w C_*(X_k),d)$
the suspension reduced simplicial complex inclusion. Regard this graded group inclusion as
an indexed chain homomorphism such that the index set has only  one element.
Then $(C_*(K;\underline{\vartheta}),d)=(\Sigma\w C_*({\cal S}(K;\underline{X},\underline{A})),d)$
and $H_*^{\XX}\!(\theta_k)=H_*^{\XX}\!(X_k,A_k)$.
\hfill$\Box$\vspace{3mm}

{\bf Example 3.10} For ${\cal S}(K;\underline{X},\underline{A})={\cal S}(K;L_1,{\cdots},L_m)$ such that all
$H_*(L_k)$ is free,
each $\iota_k\colon\w H_*(L_k)\to\w H_*(\Delta\!^{n_k})\,(=\!0)$ is an epimorphism.
By definition, $H_*^{\emptyset,\emptyset}(\Delta\!^{n_k},L_k)=0$,
$H_*^{\emptyset,\{1\}}(\Delta\!^{n_k},L_k)=\w H_{*-1}(L_k)$,
$$H_*^{\RR_m}(\underline{X},\underline{A})=H_*^{\emptyset,[m]}(\underline{X},\underline{A})=\w H_{*-1}(L_1){\otimes}{\cdots}{\otimes}\w H_{*-1}(L_m),$$
$$H_*^{\RR_m}(K)\,\widehat\otimes\,H_*^{\RR_m}(\underline{X},\underline{A})=\w H_{*-1}(K){\otimes}\w H_{*-1}(L_1){\cdots}{\otimes}\w H_{*-1}(L_m).$$
So by Theorem 3.9,
$$\w H_{*-1}({\cal S}(K;L_1,{\cdots},L_m))\cong\w H_{*-1}(K){\otimes}\w H_{*-1}(L_1){\cdots}{\otimes}\w H_{*-1}(L_m),$$
$$\w H^{*-1}({\cal S}(K;L_1,{\cdots},L_m))\cong\w H^{*-1}(K){\otimes}\w H^{*-1}(L_1){\cdots}{\otimes}\w H^{*-1}(L_m).$$

If $K$ and all $L_k$ are homology spheres ($\w H_{*}(-)\cong\Bbb Z$, so $\{\emptyset\}$
is a homology sphere but $\{\,\}$ is not),
then ${\cal S}(K;L_1,{\cdots},L_m)$ is a homology sphere. The homotopy type
of ${\cal S}(K;L_1,{\cdots},L_m)$ is discussed in \cite{AY}.
\vspace{3mm}

We have ring isomorphism $H^*({\cal S}(K;\underline{X},\underline{A}))
\cong H^*(|{\cal S}(K;\underline{X},\underline{A})|)$, where $|\cdot|$ means geometrical
realization.
So $\w H^{*}({\cal S}(K;\underline{X},\underline{A}))$ is
a ring by adding a unit to it.
This ring is not considered in this paper. \vspace{3mm}

{\bf Theorem 3.11}{\it\, For a total homology split ${\cal S}(K;\underline{X},\underline{A})$, we have
$$\begin{array}{l}
\quad H_*^{\SS_n}({\cal S}(K;\underline{X},\underline{A}))\vspace{2mm}\\
\cong H_*^{\XX_m}(K)\,\widehat\otimes\,
\big(H_*^{\XX_m;\SS_n}(\underline{X},\underline{A})\big)\vspace{2mm}\\
\cong H_*^{\XX_m}(K)\,\widehat\otimes\,
\big(H_*^{\XX\!;\SS_{n_1}}(X_1,A_1){\otimes}{\cdots}{\otimes}H_*^{\XX\!;\SS_{n_m}}(X_m,A_m)\big),
\end{array}$$
$$\begin{array}{l}
\quad H^{\,*}_{\!\SS_n}({\cal S}(K;\underline{X},\underline{A}))\vspace{2mm}\\
\cong H^{\,*}_{\!\XX_m}(K)\,\widehat\otimes\,
\big(H^{\,*}_{\!\XX_m;\SS_n}(\underline{X},\underline{A})\big)\vspace{2mm}\\
\cong H^{\,*}_{\!\XX_m}(K)\,\widehat\otimes\,
\big(H^{\,*}_{\!\XX\!;\SS_{n_1}}(X_1,A_1){\otimes}{\cdots}{\otimes}H^{\,*}_{\!\XX\!;\SS_{n_m}}(X_m,A_m)\big),  \end{array}
$$
where $\SS_n=\SS_{n_1}{\times}{\cdots}{\times}\SS_{n_m}$,
$\SS=\XX$ or $\RR$ ($\SS_{n_1}=\XX_{n_1}$, $\SS_{n_2}=\RR_{n_2}$ is possible).
For $(\hat\sigma,\hat\omega)\in\XX_m$, $(\sigma,\omega)\in\SS_n$,
$\sigma_k=\sigma{\cap}[n_k]$, $\omega_k=\omega{\cap}[n_k]$,
$$H_*^{s;\sigma_k,\omega_k}(X_k,A_k)=
\left\{\begin{array}{cl}
\Sigma{\rm coker}\,\iota_{\sigma_k,\omega_k}&{\rm if}\, s={\scriptstyle\{\{1\},\emptyset\}},\vspace{1mm}\\
\Sigma{\rm ker}\,\iota_{\sigma_k,\omega_k}&{\rm if}\, s={\scriptstyle\{\emptyset,\{1\}\}},\vspace{1mm}\\
\Sigma{\rm im}\,\iota_{\sigma_k,\omega_k}&{\rm if}\, s={\scriptstyle\{\,\emptyset,\,\emptyset\}},
\end{array}\right.$$
$$H^*_{s;\sigma_k,\omega_k}(X_k,A_k)=
\left\{\begin{array}{cl}
\Sigma{\rm ker}\,\iota_{\sigma_k,\omega_k}^\circ&{\rm if}\, s={\scriptstyle\{\{1\},\emptyset\}},\vspace{1mm}\\
\Sigma{\rm coker}\,\iota_{\sigma_k,\omega_k}^\circ&{\rm if}\, s={\scriptstyle\{\emptyset,\{1\}\}},\vspace{1mm}\\
\Sigma{\rm im}\,\iota_{\sigma_k,\omega_k}^\circ&{\rm if}\, s={\scriptstyle\{\,\emptyset,\,\emptyset\}},
\end{array}\right.$$
$$H_*^{\hat\sigma\!,\,\hat\omega;\sigma\!,\,\omega}(\underline{X},\underline{A})
=H_1{\otimes}{\cdots}{\otimes}H_m,\quad
H_k=\left\{\begin{array}{rl}
\Sigma{\rm coker}\,\iota_{\sigma_k,\omega_k}&{\rm if}\, k\in\hat\sigma,\vspace{1mm}\\
\Sigma{\rm ker}\,\iota_{\sigma_k,\omega_k}&{\rm if}\, k\in\hat\omega,\vspace{1mm}\\
\Sigma{\rm im}\,\iota_{\sigma_k,\omega_k}
&{\rm otherwise},
\end{array}\right.$$
$$H^*_{\hat\sigma\!,\,\hat\omega;\sigma\!,\,\omega}(\underline{X},\underline{A})
=H^1{\otimes}{\cdots}{\otimes}H^m,\quad
H^k=\left\{\begin{array}{rl}
\Sigma{\rm ker}\,\iota_{\sigma_k,\omega_k}^\circ&{\rm if}\, k\in\hat\sigma,\vspace{1mm}\\
\Sigma{\rm coker}\,\iota_{\sigma_k,\omega_k}^\circ&{\rm if}\, k\in\hat\omega,\vspace{1mm}\\
\Sigma{\rm im}\,\iota_{\sigma_k,\omega_k}^\circ&{\rm otherwise},
\end{array}\right.$$
where $\iota_{\sigma_k,\omega_k}$ is as in Definition 3.8
with dual $\iota_{\sigma_k,\omega_k}^\circ$ and $\Sigma$ means suspension.

If each $\iota_{\sigma_k,\omega_k}$ is an epimorphism, then all $\XX$ is replaced by $\RR$ and we have
$H_*^{\RR;\SS_{n_k}}(X_k,A_k)=H_*^{\SS_{n_k}}(A_k)$,
$H^{\,*}_{\RR;\SS_{n_k}}(X_k,A_k)=H^{\,*}_{\SS_{n_k}}(A_k)$.

If the (right) total (co)homology group is taken over a field, then the theorem holds for all
polyhedral product complexes.
\vspace{2mm}

Proof}\, A corollary of Theorem 3.7
by taking $\underline{\vartheta}=\{\vartheta_k\}_{k=1}^m$ with split inclusion
$\vartheta_k\colon T_*^{\SS_{n_k}}(A_k)\to T_*^{\SS_{n_k}}(X_k)$
the (right) total chain complex (as in Definition 4.6 in \cite{Z}) inclusion.
Then $(C_*(K;\underline{\vartheta}),d)=(T_*^{\SS_n}({\cal S}(K,\underline{X},\underline{A})),d)$
and $H_*^\XX\!(\theta_k)=H_*^{\XX\!;\SS_k}(\theta_k)=H_*^{\XX\!;\SS_k}(X_k,A_k)$.
\hfill$\Box$\vspace{3mm}

{\bf Example 3.12}  Apply Theorem 3.11 for ${\cal S}(K;\underline{X},\underline{A})={\cal S}(K;L_1,{\cdots},L_m)$,
$L_k\neq\{\,\}$ for $k=1,{\cdots},m$.
So either all $H_*^{\RR_{n_k}}(L_k)$ are free or the (co)homology is taken over a field.

For $\SS=\RR$, $\theta_k\colon H_*^{\RR_{n_k}}(L_k)\to H_*^{\RR_{n_k}}(\Delta\!^{n_k})\cong \Bbb Z$
is an epimorphism. So
$$H_*^{\RR_n}({\cal S}(K;L_1,{\cdots},L_m))\cong H_*^{\RR_m}(K)\,\widehat\otimes\,
\big(H_*^{\RR_{n_1}}(L_1){\otimes}{\cdots}{\otimes}H_*^{\RR_{n_m}}(L_m)\big),$$
$$H^*_{\RR_n}({\cal S}(K;L_1,{\cdots},L_m))\cong H^*_{\RR_m}(K)\,\widehat\otimes\,
\big(H^*_{\RR_{n_1}}(L_1){\otimes}{\cdots}{\otimes}H^*_{\RR_{n_m}}(L_m)\big).$$

For $\SS=\XX$, $\theta_k\colon H_*^{\XX_{n_k}}\!(L_k)\to H_*^{\XX_{n_k}}\!(\Delta\!^{n_k})\cong \Bbb Z$
is not an epimorphism.
To simplify notation, $H_*^{-;-}(\Delta\!^{n_k}{,}L_k)$ is abbreviated to $H_*^{-;-}$.
For $\omega_k\neq\emptyset$, $H_*^{\sigma_k,\omega_k}(\Delta\!^{n_k})=0$.
For $\sigma_k\in L_k$, $H_*^{\sigma_k,\emptyset}(L_k)\cong H_*^{\sigma_k,\emptyset}(\Delta\!^{n_k})\cong\Bbb Z$.
For $\sigma_k\notin L_k$, $H_*^{\sigma_k,\emptyset}(L_k)=0$, $H_*^{\sigma_k,\omega_k}(\Delta\!^{n_k})\cong\Bbb Z$.
So $H_*^{\emptyset,\{1\};\sigma_k,\omega_k}\cong H_*^{\sigma_k,\omega_k}(L_k)$ for $\omega_k\neq\emptyset$,
$H_*^{\emptyset,\emptyset;\sigma_k,\emptyset}\cong\Bbb Z$ for $\sigma_k\in L_k$ and
$H_*^{\{1\},\emptyset;\sigma_k,\emptyset}\cong\Bbb Z$ for $\sigma_k\notin L_k$. So
we have that for $(\sigma,\omega)\in\XX_n$,
$$H_*^{\sigma,\omega}({\cal S}(K;L_1,{\cdots},L_m))\cong H_*^{\hat\sigma,\hat\omega}(K)\otimes
\big(\otimes_{\omega_k\neq\emptyset}\,H_*^{\sigma_k,\omega_k}(L_k)\big),$$
$$H^*_{\sigma,\omega}({\cal S}(K;L_1,{\cdots},L_m))\cong H^*_{\hat\sigma,\hat\omega}(K)\otimes
\big(\otimes_{\omega_k\neq\emptyset}\,H^*_{\sigma_k,\omega_k}(L_k)\big),$$
where $\hat\sigma=\{k\,|\,\sigma_k\!\notin\! L_k,\,\omega_k=\emptyset\}$, $\hat\omega=\{k\,|\,\omega_k\!\neq\!\emptyset\}$,
$\sigma_k=\sigma{\cap}[n_k]$, $\omega_k=\omega{\cap}[n_k]$.

\section{Universal Algebra}

In this section, we compute the (right) universal (normal, etc.) algebra of total homology split polyhedral product complexes.
The (co)associativity is not required for a (co)algebra as in \cite{Z}.
\vspace{3mm}

{\bf Theorem 4.1}\, {\it Let $\vartheta\colon (C_*,d)\to(D_*,d)$ be a split inclusion
with induced homology homomorphism $\theta\colon U_*\to V_*$
such that $\vartheta$ is also a coalgebra homomorphism $\vartheta\colon (C_*,\psi_C)\to(D_*,\psi_D)$
with induced homology coalgebra homomorphism $\theta\colon (U_*,\psi_U)\to(V_*,\psi_V)$.

Then the group $C_*^{\XX}\!(\theta)$ in Definition 3.3 has a unique character coproduct
$\widehat\psi(\vartheta)$ satisfying the following three conditions.

i) $\widehat\psi(\vartheta)$ makes the following diagram ($q$, $q'$ and $\vartheta'$ as in Theorem 3.4).

\begin{center}
\setlength{\unitlength}{1mm}
\begin{picture}(79,37)
\put(0,20){${\scriptstyle C_*}$}\put(30,20){${\scriptstyle U_*}$}
\put(10,30){${\scriptstyle C_*{\otimes}C_*}$}\put(42,30){${\scriptstyle U_*{\otimes}U_*}$}
\put(0,0){${\scriptstyle D_*}$}\put(30,0){${\scriptstyle C_*^{\XX\!}}$}
\put(10,10){${\scriptstyle D_*{\otimes}D_*}$}\put(41,10){${\scriptstyle C_*^{\XX\!}{\otimes}C_*^{\XX\!}}$}
\thicklines
\put(6,21){\vector(1,0){22}}\put(12,20.5){$^{q'}$} \put(2,18){\vector(0,-1){13}}\put(2.5,12){$^\vartheta$}
\put(6,1){\vector(1,0){22}}\put(12,0.5){$^{q}$} \put(32,18){\vector(0,-1){13}}\put(32.5,12){$^{\vartheta'}$}
\thinlines
\put(23,31){\vector(1,0){15}}\put(24,30.5){$^{q'{\otimes}q'}$}
\put(15.5,28){\line(0,-1){6}}\put(15.5,20){\vector(0,-1){6}}\put(16,22){$^{\vartheta{\otimes}\vartheta}$}
\put(23,11){\line(1,0){8}}\put(33,11){\vector(1,0){6}}\put(24,10.5){$^{q{\otimes}q}$}
\put(47.5,28){\vector(0,-1){14}}\put(48,22){$^{\vartheta'{\otimes}\vartheta'}$}
\thicklines
\put(4,24){\vector(1,1){5}}\put(7,23.5){$^{\psi_C}$} \put(35,24){\vector(1,1){5}}\put(38,23){$^{\psi_U}$}
\put(4,3.5){\vector(1,1){5}}\put(7,2.5){$^{\psi_D}$} \put(35,3.5){\vector(1,1){5}}\put(38,2){$^{\widehat\psi(\vartheta)}$}
\put(57,0.5){${\scriptstyle (C_*^\XX=\,C_*^{\XX}\!(\theta))}$}
\end{picture}\end{center}
commutative except the homotopy commutative $(q{\otimes}q)\psi_D\simeq \widehat\psi(\vartheta)q$.

ii) $\widehat\psi(\vartheta)$ is independent of the choice of $\psi_C,\psi_D$ up to homotopy,
i.e., if $\psi_C,\psi_D$ are replaced by $\psi'_C,\psi'_D$ such that $\psi'_C\simeq\psi_C$, $\psi'_D\simeq\psi_D$
and we get $\widehat\psi'(\vartheta)$ for $\psi'_C$ and $\psi'_D$,
then $\widehat\psi'(\vartheta)=\widehat\psi(\vartheta)$.

iii) Denote by $\alpha={\rm coker}\,\theta$, $\beta=\Sigma\,{\rm ker}\,\theta$,
$\gamma={\rm ker}\,\theta$, $\eta={\rm im}\,\theta$.
Then $\widehat\psi(\vartheta)$ satisfies the following four conditions.

(1) $\widehat\psi(\vartheta)(\eta)\subset \eta{\otimes}\eta\oplus\gamma{\otimes}\eta\oplus \eta{\otimes}\gamma
\oplus \gamma{\otimes}\gamma$.

(2) $\widehat\psi(\vartheta)(\gamma)\subset \gamma{\otimes}\gamma\oplus \gamma{\otimes}\eta\oplus \eta{\otimes}\gamma$.

(3) $\widehat\psi(\vartheta)(\beta)\subset\big(\beta{\otimes}\gamma\oplus \beta{\otimes}\eta\oplus\eta{\otimes}\beta\big)
\oplus\big(\alpha{\otimes}\alpha\oplus \alpha{\otimes}\eta\oplus
\eta{\otimes}\alpha\oplus \eta{\otimes}\eta\big)$.

(4) $\widehat\psi(\vartheta)(\alpha)\subset \alpha{\otimes}\alpha\oplus \alpha{\otimes}\eta\oplus \eta{\otimes}\alpha
\oplus \eta{\otimes}\eta$.
\vspace{2mm}

Proof}\, Denote $N\oplus(\alpha{\oplus}\eta){\otimes}(\alpha{\oplus}\eta)=C_*^\XX{\otimes}C_*^\XX$.
By K\"{u}nneth theorem, $H_*(N)=0$ and $H_*(C_*^\XX{\otimes}C_*^\XX)=(\alpha{\oplus}\eta){\otimes}(\alpha{\oplus}\eta)$.

From the construction of $q$ and $q'$ in Theorem 3.4 we have that there is a subcomplex $F_*$ of $D_*$ such that
$H_*(F_*)=0$ and
$(D_*,d)=(F_*,d)\oplus(C_*^\XX,d)$. $q$ is the projection such that
$\vartheta({\rm ker}q')\subset F_*$.
This implies that if we take $f$ to be the composite
$$f\colon C_*^\XX\subset D_*\stackrel{\psi_D}{\longrightarrow}D_*{\otimes}D_*
\stackrel{q\otimes q}{-\!\!-\!\!\!\longrightarrow}C_*^\XX{\otimes}C_*^\XX,$$
then $f$ makes the cubic diagram of the theorem commutative except the homotopy commutative
$(q{\otimes}q)\psi_D\simeq fq$.

Construct coproduct chain homomorphism
$\psi\colon(C_*^\XX,d)\to(C_*^\XX{\otimes}C_*^\XX,d)$
and homotopy $s\colon(C_*^\XX,d)\to(\Sigma C_*^\XX{\otimes}C_*^\XX,d)$
such that $ds+sd=f{-}\psi$ and $\psi$ is independent of the choice of $f$
just as in the proof of Theorem 2.8 in \cite{Z} by regarding the graded group $\alpha,\beta,\gamma,\eta$
in that theorem as the indexed groups with the same symbol in this theorem.
\hfill$\Box$\vspace{3mm}

{\bf Definition 4.2} For the $\vartheta$ in Theorem 4.1, all the chain complexes in Theorem 3.4 are coalgebras defined as follows.

The {\it character coalgebra complex of $\vartheta$} is $(C_*^\XX\!(\theta),\widehat\psi(\vartheta))$.

If $\theta$ is an epimorphism, then $C_*^\XX(\theta)=C_*^\RR(\theta)$ and $(C_*^\RR(\theta),\widehat\psi(\vartheta))$
is called the {\it right character coalgebra complex of $\vartheta$}.

The {\it homology coalgebra of $\vartheta$} is $(H_*^{\XX}\!(\theta),\psi(\vartheta))$ with coproduct
defined as follows.

(1) $\psi(\vartheta)(x)=\widehat\psi(\vartheta)(x)$ for all $x\in \alpha\oplus\eta$.

(2) For a generator $x\in\gamma$, there is a unique generator $\overline x\in\beta$ such that $d\,\overline x= x$.
Suppose $\widehat\psi(\vartheta)(\overline x)= z+y$ with $z\in\beta{\otimes}\gamma\oplus \beta{\otimes}\eta\oplus\eta{\otimes}\beta$
and $y\in\alpha{\otimes}\alpha\oplus \alpha{\otimes}\eta\oplus\eta{\otimes}\alpha\oplus \eta{\otimes}\eta$.
Then define $\psi(\vartheta)(x)=\widehat\psi(\vartheta)(x)+y$.

If $\theta$ is an epimorphism, then $(H_*^\RR(\theta)({=}H_*^\XX\!(\theta)),\psi(\vartheta))$ is the {\it right homology coalgebra of $\vartheta$}.
The group isomorphisms $H_*^\RR(\theta)\cong V_*$ is in general not an algebra isomorphism.

For $\SS=\XX$ or $\SS$, the dual algebra $(H^{\,*}_{\RR}(\theta^\circ),\pi(\vartheta^\circ))$)
of $(H_*^\SS(\theta),\psi(\vartheta))$ is the {\it (right) cohomology algebra of $\vartheta$}.

The {\it index coalgebra complex
$(T_*^\XX\!,\widehat\psi_\vartheta)$ of $\vartheta$} is defined as follows.
Let symbols $x,x'_1,x''_2,\cdots$ be one of $\alpha,\beta,\gamma,\eta$
that are both the generators of $T_*^\XX$ and the summand groups of $C_*^\XX\!(\theta)$ in
Theorem 4.1. If the summand group $x$ satisfies $\widehat\psi(\vartheta)(x)\subset\oplus_i\,x'_i{\otimes}x''_i$
such that no summand $x'_i{\otimes}x''_i$ can be canceled,
then the generator $x$ satisfies $\widehat\psi_\vartheta(x)=\Sigma_i\,x'_i{\otimes}x''_i$.

If $\theta$ is an epimorphism, then $(T_*^\RR,\psi_{\vartheta})$ is
the {\it right index coalgebra complex of $\vartheta$}.
\vspace{3mm}

{\bf Theorem 4.3} {\it For the coalgebras in Definition 4.2,
all the chain homomorphisms in Theorem 3.4 induce cohomology algebra isomorphisms.
\vspace{2mm}

Proof}\, $q$ and $q'$ induce isomorphisms by Theorem 4.1. The proof that $\phi$
is a coalgebra isomorphism is the same as that of Theorem 6.4 in \cite{Z}.
\hfill$\Box$\vspace{3mm}

{\bf Definition 4.4} Let $K$, $\underline{\vartheta}$, $\underline{\theta}$ be
as in Definition 3.5 and Definition 3.6
such that each $\vartheta_k\colon ((C_k)_*,\psi_{C_k})\to((D_k)_*,\psi_{D_k})$
and $\theta_k\colon ((U_k)_*,\psi_{U_k})\to((V_k)_*,\psi_{V_k})$
satisfy the condition of Theorem 4.1. $\SS=\XX$ or $\RR$.
Then all the chain complexes in Theorem 3.7 are coalgebras defined as follows.

The subgroup $C_*(K;\underline{\vartheta})$ of $(D_1)_*{\otimes}{\cdots}{\otimes}(D_m)_*$ is a subcoalgebra and
is called the {\it (right) coalgebra complex of $(K;\underline{\vartheta})$}. Its cohomology algebra is denoted by
$(H^*(C^*(K;\underline{\vartheta}^\circ)),\cup_{(K;\underline{\vartheta}^\circ)})$.

Define $(C_*^{\SS_m}(\underline{\vartheta}),\widehat\psi(\underline{\vartheta}))=
(C_*^{\SS}\!(\vartheta_1)\otimes{\cdots}\otimes C_*^{\SS}\!(\vartheta_m),\widehat\psi(\vartheta_1){\otimes}{\cdots}{\otimes}
\widehat\psi(\vartheta_m))$.
Then the subgroup $C_*^{\SS_m}(K;\underline{\vartheta})$ of $C_*^{\SS_m}(\underline{\vartheta})$ is a subcoalgebra and
is called the {\it (right) character coalgebra complex of $(K;\underline{\vartheta})$}.

The {\it (right) homology coalgebra of $\underline{\vartheta}$} is
$(H_*^{\SS_m}(\underline{\theta}),\psi(\underline{\vartheta}))=
(H_*^{\SS}\!(\theta_1)\otimes{\cdots}\otimes H_*^{\SS}\!(\theta_m),\psi(\vartheta_1){\otimes}{\cdots}{\otimes}
\psi(\vartheta_m))$. The {\it (right) cohomology algebra of $\underline{\vartheta}$} is the dual algebra
$(H^*_{\SS_m}(\underline{\theta}^\circ),\pi(\underline{\vartheta}^\circ))=
(H^*_{\SS}(\theta_1^\circ)\otimes{\cdots}\otimes H^*_{\SS}(\theta_m^\circ),\pi(\vartheta_1^\circ){\otimes}{\cdots}{\otimes}
\pi(\vartheta_m^\circ))$.

Define $(T_*^{\SS_m},\widehat\psi_{\underline{\vartheta}})=
(T_*^{\SS}\otimes{\cdots}\otimes T_*^{\SS},\widehat\psi_{\vartheta_1}{\otimes}{\cdots}{\otimes}
\widehat\psi_{\vartheta_m})$.
Then the subgroup $T_*^{\SS_m}(K)$ of $T_*^{\SS_m}$ (as defined in Definition 4.5 in \cite{Z}) is a subcoalgebra
and is called the {\it (right) index coalgebra complex of $K$ induced by $\underline{\vartheta}$}.
Its cohomology algebra is called the {\it (right) index cohomology algebra of $K$ induced by $\underline{\vartheta}$} and
is denoted by $(H^*_{\SS_m}(K),\cup_{(K;\underline{\vartheta}^\circ)})$.
\vspace{3mm}

{\bf Theorem 4.5} {\it For the coalgebras in Definition 4.4,
all the chain homomorphisms in Theorem 3.7 induce cohomology algebra isomorphisms.
So we have cohomology algebra isomorphism
$$(H^*(C^*(K;\underline{\theta}^\circ)),\cup_{(K;\underline{\vartheta}^\circ)})
\cong (H^*_{\XX_m}(K)\,\widehat\otimes\, H^*_{\XX_m}(\underline{\theta}^\circ),\cup_{(K;\underline{\vartheta}^\circ)}\,\widehat\otimes\,
\pi(\underline{\vartheta}^\circ)).$$

Proof}\, The $\varphi_\sigma$ and $\phi_\sigma$ in the proof of Theorem 3.7 as a tensor product induce isomorphisms by
Theorem 4.3.
So $\varphi_{(K,\underline{\vartheta})}=+_{\sigma\in K}\,\varphi_\sigma$ and $\phi_{(K,\underline{\theta})}=+_{\sigma\in K}\,\phi_\sigma$
also induce isomorphisms.
\hfill$\Box$\vspace{3mm}

{\bf Example 4.6} Let everything be as in Theorem 2.9.
We compute the cohomology algebra of ${\cal Z}(K;\underline{Y},\underline{B})$ over a field.
By Theorem 6.9 and 7.11 in \cite{Z}, we have
$$\begin{array}{c}
(H^*(Y_k),\cup)\cong\big(H^*_{\SS_{n_k}}(X_k)\,\widehat\otimes\,H^*_{\SS_{n_k}}(\underline{U_k},\underline{C_k}),
\cup_{X_k}\,\widehat\otimes\,\pi_{(\underline{U_k},\underline{C_k})}\big),\vspace{2mm}\\
(H^*(B_k),\cup)\cong\big(H^*_{\SS_{n_k}}(A_k)\,\widehat\otimes\,H^*_{\SS_{n_k}}(\underline{U_k},\underline{C_k}),
\cup_{A_k}\,\widehat\otimes\,\pi_{(\underline{U_k},\underline{C_k})}\big),
\end{array}$$
where $\SS_{n_k}=\XX_{n_k}$ or $\RR_{n_k}$ (if $H_*(C_i)\to H_*(U_i)$ induced by inclusion is an epimorphism for $i\in[n_k]$)
and $\cup_{X_k}$ and $\cup_{A_k}$ may be a (right) universal (or (right) normal, etc.) product appearing in the theorems.

Take $\vartheta_k\colon(T_*^{\SS_{n_k}}(A_k),\psi_k)\to(T_*^{\SS_{n_k}}(X_k),\psi_k)$ to be the (right) total
chain complex inclusion and apply Theorem 4.5 for $\underline{\vartheta}=\{\vartheta_k\}_{k=1}^m$.
Then we have
$$\,\,(H^*({\cal Z}(K;\underline{Y},\underline{B}),\cup)\cong\big(H^*_{\SS_n}({\cal S}(K;\underline{X},\underline{A}))\,\widehat\otimes\,
H^*_{\SS_n}(\underline{U},\underline{C}),
\cup\,\widehat\otimes\,\pi_{(\underline{U},\underline{C})}\big),$$
$$(H^*_{\SS_n}({\cal S}(K;\underline{X},\underline{A})),\cup)\cong
\big(H^*_{\XX_m}(K)\,\widehat\otimes\,H^*_{\XX_m;\SS_n}(\underline{X},\underline{A}),
\cup_{(K;\underline{\vartheta}^\circ)}\,\widehat\otimes\,\pi(\underline{\vartheta}^\circ)\big),$$
where $\SS_n=\SS_{n_1}{\times}{\cdots}{\times}\SS_{n_m}$ ($\SS_{n_1}=\XX_{n_1}$, $\SS_{n_2}=\RR_{n_2}$ is possible).

Specifically, take each $(X_k,A_k)=(\Delta\!^{n_k},L_k)$ and suppose $H_*(C_i){\to}H_*(U_i)$
induced by inclusion is an epimorphism for all $i$.
Then the coalgebra homomorphism $\theta_k\colon H_*^{\RR_m}(L_k)\to H_*^{\RR_m}(\Delta\!^{n_k})\cong\Bbb Z$
is an epimorphism. By Theorem 4.1 for $\vartheta=\vartheta_k$,
$\eta\cong H^*_{\emptyset,\emptyset}(L_k)\cong H^*_{\emptyset,\emptyset}(\Delta\!^{n_k})\cong\Bbb Z$,
$\gamma\cong \oplus_{\omega_k\neq\emptyset}H^*_{\emptyset,\omega_k}(L_k)$, $\alpha=0$.
Since $\eta{\otimes}\eta$ is non-zero only at degree $0$ and $\beta$ is zero at degree $0$,
we have that $\eta{\otimes}\eta$ can not be a summand of $\widehat\psi(\vartheta_k)(\beta)$. So

$\widehat\psi(\vartheta_k)(\eta)\subset\eta{\otimes}\eta$.

$\widehat\psi(\vartheta_k)(\gamma)\subset\gamma{\otimes}\gamma\oplus\gamma{\otimes}\eta\oplus\eta{\otimes}\gamma$.

$\widehat\psi(\vartheta_k)(\beta)\subset\beta{\otimes}\gamma\oplus\beta{\otimes}\eta\oplus\eta{\otimes}\beta$.

So  $\widehat\psi_{\vartheta_k}$ is the right strictly normal coproduct in Definition 7.4 in \cite{Z}
and the group isomorphism $H_*^{\RR;\RR_{n_k}}(\Delta\!^{n_k},L_k)\cong H_*^{\RR_{n_k}}(L_k)$ is an algebra isomorphism.
So we have cohomology algebra isomorphism
$$H^*_{\RR_n}({\cal S}(K;L_1,{\cdots},L_m))
\cong H^*_{\RR_m}(K){\otimes}_{\RR_m}
\big(H^*_{\RR_{n_1}}(L_1){\otimes}{\cdots}{\otimes}H^*_{\RR_{n_m}}(L_m)\big),$$
where $H^*_{\RR_m}(K)$ is the right strictly normal algebra of $K$.
\vspace{3mm}

\section{Duality Isomorphism}\vspace{3mm}

In this section, we compute the Alexander duality isomorphism
on some special type of polyhedral product spaces.\vspace{3mm}

{\bf Theorem 5.1} {\it Let $(\underline{X},\underline{A})=\{(X_k,A_k)\}_{k=1}^m$ be a sequence of topological
pairs satisfying the following conditions.

1) Each homology group homomorphism $i_k\colon H_*(A_k)\to H_*(X_k)$
induced by inclusion is a split homomorphism.

2) Each $X_k$ is a closed orientable manifold of dimension $r_k$.

3) Each $A_k$ is a proper compact polyhedron subspace of $X_k$.

Let $(\underline{X},\underline{A}^c)=\{(X_k,A_k^c)\}_{k=1}^m$ with $A_k^c=X_k{\setminus}A_k$.
Then for all $(\sigma\!,\omega)\in\XX_m$, there are duality isomorphisms (\,$r=r_1{+}{\cdots}{+}r_m$\,)
$$\gamma_{\sigma\!,\,\omega}\colon H_*^{\sigma\!,\,\omega}(\underline{X},\underline{A})\to H^{r-|\omega|-*}_{\tilde\sigma\!,\,\omega}(\underline{X},\underline{A}^c),$$
$$\gamma^*_{\sigma\!,\,\omega}\colon H^*_{\sigma\!,\,\omega}(\underline{X},\underline{A})\to H_{r-|\omega|-*}^{\tilde\sigma\!,\,\omega}(\underline{X},\underline{A}^c),$$
where $\w\sigma=[m]{\setminus}(\sigma{\cup}\omega)$,
$H_*^{\sigma\!,\,\omega}(-)$ and $H^*_{\sigma\!,\,\omega}(-)$ are as in Theorem 3.9.

If the (co)homology is taken over a field, then the conclusion holds for
$(\underline{X},\underline{A})$ satisfying the following conditions.

1) Each $X_k$ is a closed manifold of dimension $r_k$ orientable with respect to the homology theory over the field.

2) Each $A_k$ is a proper compact polyhedron subspace of $X_k$.
}\vspace{2mm}

{\it Proof}\, We have the following commutative diagram of exact sequences
\[\begin{array}{cccccccccc}
\cdots\longrightarrow \hspace{-1.5mm}&\hspace{-1.5mm} {\scriptstyle H_n(A_k)} \hspace{-1.5mm}&\hspace{-1.5mm} \stackrel{i_k}{\longrightarrow} \hspace{-1.5mm}&\hspace{-1.5mm} {\scriptstyle H_n(X_k)}
&\hspace{-1.5mm} \stackrel{j_k}{\longrightarrow} \hspace{-1.5mm}&\hspace{-1.5mm} {\scriptstyle H_n(X_k,A_k)} \hspace{-1.5mm}&\hspace{-1.5mm}\stackrel{\partial_k}{\longrightarrow} \hspace{-1.5mm}&\hspace{-1.5mm}
{\scriptstyle H_{n-1}(A_k)}\hspace{-1.5mm}&\hspace{-1.5mm} \longrightarrow \cdots\,\vspace{1mm}\\
&^{\alpha_k}\downarrow\quad&&^{\gamma_k}\downarrow\quad&&^{\beta_k}\downarrow\quad&&^{\alpha_k}\downarrow\quad\quad&&\\
\cdots\longrightarrow \hspace{-1.5mm}&\hspace{-1.5mm} {\scriptstyle H^{r_k-n}(X_k,A_k^c)} \hspace{-1.5mm}&\hspace{-1.5mm} \stackrel{q^*_k}{\longrightarrow} \hspace{-1.5mm}&\hspace{-1.5mm}{\scriptstyle H^{r_k-n}(X_k)}
\hspace{-1.5mm}&\hspace{-1.5mm} \stackrel{p^*_k}{\longrightarrow} \hspace{-1.5mm}&\hspace{-1.5mm}
{\scriptstyle H^{r_k-n}(A_k^c)} \hspace{-1.5mm}&\hspace{-1.5mm}\stackrel{\partial^*_k}{\longrightarrow} \hspace{-1.5mm}&\hspace{-1.5mm}
{\scriptstyle H^{r_k-n+1}(X_k,A^c_k)} \hspace{-1.5mm}&\hspace{-1.5mm}\longrightarrow\cdots,
  \end{array}\]
where $\alpha_k,\beta_k$ are the Alexander duality isomorphisms
and $\gamma_k$ is the Poncar\'{e} duality isomorphism.
So we have the following group isomorphisms
$$\begin{array}{rl}
(\partial^*_k)^{-1}\alpha_k\colon&\,\,\,{\rm ker}\,i_k\,\,\,\stackrel{\cong}
{\longrightarrow}\,\,{\rm coker}\,p_k^*,\vspace{1mm}\\
\gamma_k\colon&\,\,\,{\rm im}\,i_k\,\,\,\,\stackrel{\cong}{\longrightarrow}
\,\,\,{\rm ker}\,p_k^*,\vspace{1mm}\\
p^*_k\gamma_k\colon&{\rm coker}\,i_k\stackrel{\cong}{\longrightarrow}
\,\,\,\,{\rm im}\,p_k^*.
\end{array}$$
Define $\theta_k\colon H_*^\XX(X_k,A_k)\to H^*_\XX(X_k,A_k^c)$
to be the direct sum of the above three isomorphisms.
Then $\theta_1{\otimes}{\cdots}{\otimes}\theta_m=\oplus_{(\sigma\!,\omega)\in\XX_m}\,\gamma_{\sigma\!,\omega}$.
\hfill $\Box$\vspace{3mm}

{\bf Theorem 5.2} {\it Let $K$ and $K^\circ$ be the dual of each other relative to $[m]$.
Then  for all $(\sigma\!,\omega)\in\XX_m$, $\omega\neq\emptyset$,
there are duality isomorphisms
$$\gamma_{K,\sigma\!,\,\omega}\colon H_*^{\sigma\!,\,\omega}(K)=\w H_{*-1}(K_{\sigma\!,\omega})\to
H^{|\omega|-*-1}_{\tilde\sigma\!,\,\omega}(K^\circ)=\w H^{|\omega|-*-2}((K^\circ)_{\tilde\sigma\!,\omega}),$$
$$\gamma^*_{K,\sigma\!,\,\omega}\colon H^*_{\sigma\!,\,\omega}(K)=\w H^{*-1}(K_{\sigma\!,\omega})
\to H_{|\omega|-*-1}^{\tilde\sigma\!,\,\omega}(K^\circ)=\w H_{|\omega|-*-2}((K^\circ)_{\tilde\sigma\!,\omega}),$$
where $\w\sigma=[m]{\setminus}(\sigma{\cup}\omega)$, $|\omega|$ is the cardinality of $\omega$.
}\vspace{2mm}

{\it Proof}\, Let $(C_*(\Delta\!^{\omega},K_{\sigma\!,\,\omega}),d)$ be the relative simplicial chain complex.
Since $\w H_*(\Delta\!^{\omega})=0$, we have a boundary isomorphism\vspace{1mm}\\
\hspace*{30mm}$\partial\colon H_*(\Delta\!^\omega,K_{\sigma\!,\,\omega})\stackrel{\cong}{\longrightarrow}
\w H_{*-1}(K_{\sigma\!,\,\omega})=H_*^{\sigma\!,\,\omega}(K)$.\vspace{1mm}\\
$C_*(\Delta\!^{\omega},K_{\sigma\!,\,\omega})$ has a set of generators consisting of all non-simplices of $K_{\sigma\!,\,\omega}$, i.e., $K^c_{\sigma\!,\,\omega}=\{\eta\subset\omega\,|\,\eta\not\in K_{\sigma\!,\,\omega}\}$ is a set of generators of $C_*(\Delta\!^{\omega},K_{\sigma\!,\,\omega})$.
So we may denote $(C_*(\Delta\!^\omega,K_{\sigma\!,\,\omega}),d)$ by $(C_*(K^c_{\sigma\!,\,\omega}),d)$, where
$\eta\in K^c_{\sigma\!,\,\omega}$ has degree $|\eta|{-}1$ with $|\eta|$ the cardinality of $\eta$.
The correspondence $\eta\to \omega{\setminus}\eta$ for all $\eta\in K^c_{\sigma\!,\,\omega}$ induces a dual complex isomorphism\vspace{1mm}\\
\hspace*{32mm}$\psi\colon(C_*(K^c_{\sigma\!,\,\omega}),d)\to
(\w C^{*}((K_{\sigma\!,\,\omega})^\circ),\delta)$.\vspace{1mm}\\
Since $(K_{\sigma\!,\,\omega})^\circ=(K^\circ)_{\tilde\sigma\!,\,\omega}$, we have induced homology group isomorphism
$\bar\psi\colon H_*(\Delta\!^\omega,K_{\sigma\!,\,\omega})\to H^{|\omega|-*-1}_{\tilde\sigma\!,\,\omega}(K^\circ)$.
Define $\gamma_{K,\sigma\!,\,\omega}=\bar\psi\partial^{-1}$.
\hfill $\Box$\vspace{3mm}

Notice that for $\sigma\in K$, $[m]{\setminus}\sigma$ may not be a simplex of $K^\circ$.
In this case, there is no isomorphism from $H_*^{\sigma,\emptyset}(K)={\Bbb Z}$
to $H^*_{\tilde\sigma,\emptyset}(K^\circ)=0$.\vspace{3mm}

{\bf Example 5.3} For the ${\cal S}(K;L_1,{\cdots},L_m)$ and index sets $\sigma,\omega,\hat\sigma,\hat\omega,\sigma_k,\omega_k$ in Example 3.12,
$\gamma_{{\cal S}(K;L_1,{\cdots},L_m),\,\sigma,\,\omega}
=\gamma_{K,\hat\sigma\!,\,\hat\omega}\otimes(\otimes_{\omega_k\neq\emptyset}\,\gamma_{L_k,\sigma_k,\omega_k})$.
\vspace{3mm}

{\bf Definition 5.4} For homology split $M={\cal Z}(K;\underline{X},\underline{A})$, let
$i\colon H_*(M)\to H_*(\w X)$ and $i^*\colon H^*(\w X)\to H^*(M)$
be the singular (co)homology homomorphism induced by the inclusion map from $M$ to
$\w X=X_1{\times}{\cdots}{\times}X_m$.
From the long exact exact sequences
$$\begin{array}{l}
\scriptstyle{\cdots\longrightarrow\,\, H_n(M)\,\,\stackrel{i}{\longrightarrow}\,\,H_n(\w X)
\,\,\stackrel{j}{\longrightarrow}\,\, H_n(\w X,M)\,\,\stackrel{\partial}{\longrightarrow}\,\, H_{n-1}(M)
\,\,\longrightarrow\cdots}\\\\
\scriptstyle{\cdots\longrightarrow\,\, H^{n-1}(M)\,\,\stackrel{\partial^*}{\longrightarrow}\,\, H^n(\w X,M)
\,\,\stackrel{j^*}{\longrightarrow}\,\, H^n(\w X)\,\,\stackrel{i^*}{\longrightarrow}\,\, H^{n}(M)\,\, \longrightarrow\cdots}
  \end{array}$$
we define
$$\hat H_*(M)\!=\!{\rm coim}\,i,\,\overline H_*(M)\!=\!{\rm ker}\,i,\,
\hat H_*(\w X{,}M)\!=\!{\rm im}\,j,\,\overline H_*(\w X{,}M)\!=\!{\rm coker}\,j,$$
$$\hat H^*(M)\!=\!{\rm im}\,i^*,\,\overline H^*(M)\!=\!{\rm coker}\,i^*,\,
\hat H^*(\w X{,}M)\!=\!{\rm coim}\,j^*,\,H^*(\w X{,}M)\!=\!{\rm ker}\,j^*.\vspace{3mm}$$

{\bf Theorem 5.5} {\it For a homology split space $M={\cal Z}(K;\underline{X},\underline{A})$,
we have the following group decompositions
$$\begin{array}{cc}
H_*(M)=\hat H_*(M){\oplus}\overline H_*(M),&
H_*(\w X,M)=\hat H_*(\w X,M){\oplus}\overline H_*(\w X,M),\vspace{2mm}\\
H^*(M)=\hat H^*(M){\oplus}\overline H^*(M),&
H^*(\w X,M)=\hat H^*(\w X,M){\oplus}\overline H^*(\w X,M)
\end{array}$$
and direct sum group decompositions
$$\begin{array}{c}
\overline H_{*+1}(\w X,M)\cong \overline H_*(M)
\cong \oplus_{(\sigma\!,\,\omega)\in\overline\XX_m}\,
H_*^{\sigma\!,\,\omega}(K)\otimes H_*^{\sigma\!,\,\omega}(\underline{X},\underline{A}),\vspace{2mm}\\
\overline H^{*+1}(\w X,M)\cong \overline H^*(M)
\cong \oplus_{(\sigma\!,\,\omega)\in\overline\XX_m}\,
H^*_{\sigma\!,\,\omega}(K)\otimes H^*_{\sigma\!,\,\omega}(\underline{X},\underline{A}),\vspace{2mm}\\
\hat H_*(M)\cong\oplus_{\sigma\in K} H_*^{\sigma,\emptyset}(\underline{X},\underline{A}),\,\,
\hat H_*(\w X,M)\cong\oplus_{\sigma\notin K} H_*^{\sigma,\emptyset}(\underline{X},\underline{A}),\vspace{2mm}\\
\hat H^*(M)\cong\oplus_{\sigma\in K} H^*_{\sigma,\emptyset}(\underline{X},\underline{A}),
\,\,\hat H^*(\w X,M)\cong\oplus_{\sigma\notin K} H^*_{\sigma,\emptyset}(\underline{X},\underline{A}),
\end{array}$$
where $\overline\XX_m=\{(\sigma\!,\omega)\in\XX_m\,|\,\omega\neq\emptyset\}$.

The conclusion holds for all polyhedral product spaces if the (co)homology group is taken over a field.
}\vspace{3mm}

{\it Proof}\, By definition, $i=\oplus_{(\sigma\!,\,\omega)\in\XX_m}\,i_{\sigma\!,\,\omega}$ with\vspace{2mm}\\
\hspace*{15mm}$i_{\sigma\!,\,\omega}\colon H_*^{\sigma\!,\,\omega}(K){\otimes}H_*^{\sigma\!,\,\omega}(\underline{X},\underline{A})
\stackrel{i{\otimes}1}{-\!\!\!\longrightarrow}
H_*^{\sigma\!,\,\omega}(\Delta\!^m){\otimes}H_*^{\sigma\!,\,\omega}(\underline{X},\underline{A})$,\vspace{1mm}\\
where $i$ is induced by inclusion and $1$ is the identity.
$H_*^{\sigma\!,\,\omega}(\Delta\!^m)=0$ if $\omega\neq\emptyset$,
$H_*^{\sigma,\emptyset}(K)={\Bbb Z}$ if $\sigma\in L$
and $H_*^{\sigma,\emptyset}(K)=0$ if $\sigma\notin L$.
So
$$\hat H_*(M)=\oplus_{\sigma\in K}
H_*^{\sigma,\emptyset}(K){\otimes}H_*^{\sigma,\emptyset}(\underline{X},\underline{A})
\cong\oplus_{\sigma\in K} H_*^{\sigma,\emptyset}(\underline{X},\underline{A})$$
$$\overline H_*(M)
=\oplus_{(\sigma\!,\,\omega)\in\overline\XX_m}\,
H_*^{\sigma\!,\,\omega}(K)\otimes H_*^{\sigma\!,\,\omega}(\underline{X},\underline{A}).$$

The relative group case is similar.
\hfill $\Box$\vspace{3mm}

{\bf Theorem 5.6} {\it For the space $M={\cal Z}(K;\underline{X},\underline{A})$ such that $(\underline{X},\underline{A})$
satisfies the condition of Theorem 5.1,
the Alexander duality isomorphisms
$$\alpha\colon H_*(M)\to H^{r-*}(\w X,M^c),\,\quad
\alpha^*\colon H^*(M)\to H_{r-*}(\w X,M^c)$$
satisfy $\,\alpha=\hat\alpha\oplus\overline\alpha$,\,\,$\alpha^*=\hat\alpha^*\oplus\overline\alpha^*$,
where
$$\hat\alpha\colon \hat H_*(M)\to \hat H^{r-*}(\w X,M^c),\,\quad
\overline\alpha\colon \overline H_*(M)\to \overline H^{r-*}(\w X,M^c)\cong\overline H^{r-*-1}(M^c),$$
$$\hat\alpha^*\colon \hat H^*(M)\to \hat H_{r-*}(\w X,M^c),\,\quad
\overline\alpha^*\colon \overline H^*(M)\to \overline H_{r-*}(\w X,M^c)\cong\overline H_{r-*-1}(M^c)\,$$
are as follows. Identify all the above groups with the direct sum groups in Theorem 5.5.
Then
$$\hat\alpha\,\,=\oplus_{\sigma\in K}\,\gamma_{\sigma,\emptyset},\,\quad
\overline\alpha\,\,=\oplus_{(\sigma\!,\omega)\in\overline\XX_m}\,\,\gamma_{K,\sigma,\omega}{\otimes}\gamma_{\sigma,\omega},$$
$$\hat\alpha^*=\oplus_{\sigma\in K}\,\gamma_{\sigma,\emptyset}^*,\,\quad
\overline\alpha^*\,\,=\oplus_{(\sigma\!,\omega)\in\overline\XX_m}\,\,\gamma_{K,\sigma,\omega}^*{\otimes}\gamma_{\sigma,\omega}^*,$$
where $\gamma_{\sigma\!,\,\omega}$, $\gamma_{\sigma\!,\,\omega}^*$ are as in Theorem 5.1
and $\gamma_{K,\sigma\!,\,\omega}$, $\gamma_{K,\sigma\!,\,\omega}^*$ are as in Theorem 5.2.
}\vspace{2mm}

{\it Proof}\, Denote by $\alpha=\alpha_M$, $\hat\alpha=\hat\alpha_M$,
$\overline\alpha=\overline\alpha_M$.
Then for $M={\cal Z}(K;\underline{X},\underline{A})$ and $N={\cal Z}(L;\underline{X},\underline{A})$,
we have the following commutative diagrams of exact sequences
\[\,\,\,\,\begin{array}{ccccccc}
{\scriptstyle\cdots\,\,\longrightarrow}\!\!&\!\!{\scriptstyle H_k(M{\cap}N)}\!\!&\!\!{\scriptstyle\longrightarrow}
\!\!&\!\!{\scriptstyle H_k(M){\oplus}H_k(N)}
\!\!&\!\!{\scriptstyle\longrightarrow}\!\!&\!\!{\scriptstyle H_k(M{\cup}N)}\!\!&\!\!{\scriptstyle\longrightarrow\,\,\cdots}\vspace{1mm}\\
\!\!&\!\!^{\alpha_{M\cap N}}\downarrow\quad\quad\!\!&\!\!\!\!&\!\!^{\alpha_M\oplus\alpha_N}\downarrow\quad\quad\quad\!\!&\!\!\!\!&\!\!^{\alpha_{M\cup N}}\downarrow\quad\quad\!\!&\!\!\\
{\scriptstyle\cdots\,\,\longrightarrow}\!\!&\!\!{\scriptstyle H^{r-k}(\w X,(M{\cap}N)^c)}\!\!&\!\!{\scriptstyle\longrightarrow}\!\!&\!\!
{\scriptstyle H^{r-k}(\w X,M^c){\oplus}H^{r-k}(\w X,N^c)}
\!\!&\!\!{\scriptstyle\longrightarrow}\!\!&\!\!{\scriptstyle H^{r-k}(\w X,(M{\cup}N)^c)}\!\!&\!\!{\scriptstyle\longrightarrow\,\,\cdots}
  \end{array}\quad\,\,\,\,(1)\vspace{1mm}\]
\[\begin{array}{ccccccc}
{\scriptstyle0\quad\longrightarrow}\!\!&\!\!{\scriptstyle \hat H_k(M{\cap}N)}\!\!&\!\!{\scriptstyle\longrightarrow}\!\!&\!\!
{\scriptstyle\hat H_k(M){\oplus}\hat H_k(N)}
\!\!&\!\!{\scriptstyle\longrightarrow}\!\!&\!\!
{\scriptstyle\hat H_k(M{\cup}N;\underline{X},\underline{A})}\!\!&\!\!{\scriptstyle\longrightarrow\quad 0}
\vspace{1mm}\\
\!\!&\!\!^{\hat\alpha_{M\cap N}}\downarrow\quad\quad\!\!&\!\!\!\!&\!\!^{\hat\alpha_M\oplus\hat\alpha_N}\downarrow\quad\quad\quad\!\!&\!\!\!\!&\!\!^{\hat\alpha_{M\cup N}}\downarrow\quad\quad\!\!&\!\!\\
{\scriptstyle0\quad\longrightarrow}\!\!&\!\!{\scriptstyle\hat H^{r-k}(\w X,(M{\cap}N)^c)}\!\!&\!\!{\scriptstyle\longrightarrow}\!\!&\!\!
{\scriptstyle\hat H^{r-k}(\w X,M^c){\oplus}\hat H^{r-k}(\w X,N^c)}
\!\!&\!\!{\scriptstyle\longrightarrow}\!\!&\!\!{\scriptstyle\hat H^{r-k}(\w X,(M{\cup}N)^c)}\!\!&\!\!{\scriptstyle\longrightarrow\quad 0}
  \end{array}\quad(2)\vspace{1mm}\]
For $(\sigma,\omega)\in\overline\XX_m$,
$A=H_l^{\sigma\!,\,\omega}(\underline{X},\underline{A})$,
$B=H^{r-|\omega|-l}_{\sigma\!,\,\omega}(\underline{X},\underline{A}^c)$,
$\gamma_1=\gamma_{K{\cap}L,\sigma\!,\omega}$, $\gamma_2=\gamma_{K,\sigma\!,\omega}{\oplus}\gamma_{L,\sigma\!,\omega}$,
$\gamma_3=\gamma_{K\cup L,\sigma\!,\omega}$, we have the commutative diagram
\[\begin{array}{ccccccc}
{\scriptstyle\cdots\,\,\longrightarrow}\!\!&\!\!{\scriptstyle H_k^{\sigma\!,\omega}(K{\cap}L){\otimes}A}\!\!&\!\!{\scriptstyle\longrightarrow}\!\!&\!\!
{\scriptstyle (H_k^{\sigma\!,\omega}(K){\oplus} H_k^{\sigma\!,\omega}(L)){\otimes}A}\!\!&\!\!
{\scriptstyle\longrightarrow}\!\!&\!\!{\scriptstyle H_k^{\sigma\!,\omega}(K{\cup}L){\otimes}A}\!\!&\!\!{\scriptstyle\longrightarrow\,\,\cdots}\vspace{1mm}\\
\!\!&\!\!^{\gamma_1\otimes\gamma_{\sigma,\omega}}
\downarrow\quad\quad\!\!&\!\!\!\!&\!\!^{\gamma_2\otimes\gamma_{\sigma,\omega}}
\downarrow\quad\quad\quad\!\!&\!\!\!\!&\!\!
^{\gamma_{3}{\otimes}\gamma_{\sigma,\omega}}\downarrow\quad\quad\!\!&\!\!\\
{\scriptstyle\cdots\,\,\longrightarrow}\!\!&\!\!
{\scriptstyle H^{|\omega|-k-1}_{\tilde\sigma\!,\omega}((K\cap L)^\circ){\otimes}B}\!\!&\!\!
{\scriptstyle\longrightarrow}\!\!&\!\!{\scriptstyle (H^{|\omega|-k-1}_{\tilde\sigma\!,\omega}(K^\circ)
{\oplus}H^{|\omega|-k-1}_{\tilde\sigma\!,\omega}(L^\circ)){\otimes}B}\!\!&\!\!{\scriptstyle\longrightarrow}\!\!&\!\!
{\scriptstyle H^{|\omega|-k-1}_{\tilde\sigma\!,\omega}((K\cup L)^*){\otimes}B}\!\!&\!\!
{\scriptstyle\longrightarrow\,\,\cdots}
  \end{array}\]
The direct sum of all the above diagrams is the following diagram.
\[\begin{array}{ccccccc}
{\scriptstyle\cdots\,\,\longrightarrow}\!\!&\!\!{\scriptstyle\overline H_k(M{\cap}N)}\!\!&\!\!{\scriptstyle\longrightarrow}\!\!&\!\!
{\scriptstyle\overline H_k(M){\oplus}\overline H_k(N)}\!\!&\!\!{\scriptstyle\longrightarrow}\!\!&\!\!
{\scriptstyle\overline H_k(M{\cup}N)}\!\!&\!\!{\scriptstyle\longrightarrow\,\,\cdots}\vspace{1mm}\\
\!\!&\!\!^{\overline\alpha_{M\cap N}}\downarrow\quad\quad\!\!&\!\!\!\!&\!\!^{\overline\alpha_M\oplus\overline\alpha_N}\downarrow\quad\quad\quad\!\!&\!\!\!\!&\!\!^{\overline\alpha_{M\cup N}}\downarrow\quad\quad\!\!&\!\!\\
{\scriptstyle\cdots\,\,\longrightarrow}\!\!&\!\!{\scriptstyle\overline H^{r-k}(\w X,(M{\cap}N)^c)}\!\!&\!\!
{\scriptstyle\longrightarrow}\!\!&\!\!{\scriptstyle\overline H^{r-k}(\w X,M^c){\oplus}\overline H^{r-k}(\w X,N^c)}
\!\!&\!\!{\scriptstyle\longrightarrow}\!\!&\!\!{\scriptstyle\overline H^{r-k}(\w X,(M{\cup}N)^c)}\!\!&\!\!
{\scriptstyle\longrightarrow\,\,\cdots}
  \end{array}\quad(3)\]

(1), (2) and (3) imply that
if the theorem holds for $M$ and $N$ and $M{\cap}N$, then it holds for $M{\cup}N$.
So by induction on the number of maximal simplices of $K$,
we only need prove the theorem for the special case that $K$ has only one maximal simplex.

Now we prove the theorem for $M={\cal Z}(\Delta\!^S;\underline{X},\underline{A})$ with $S\subset[m]$.
Then
$$M=Y_1{\times}{\cdots}{\times}Y_m,\quad Y_k=\left\{\begin{array}{cc}
X_k & {\rm if}\,\,k\in S, \\
A_k & {\rm if}\,\,k\notin S.
\end{array}\right.$$
So $(\w X,M^c)=(X_1,Y_1^c){\times}{\cdots}{\times}(X_m,Y_m^c)$.

By identifying ${\rm coim}\,q^*_k$ and $\Sigma^{-1}{\rm im}\,\partial^*_k$ respectively with ${\rm ker}\,p^*_k$
and ${\rm coker}\,p^*_k$ in the following commutative diagram
\[\begin{array}{cccccccccc}
\cdots\longrightarrow \hspace{-1.5mm}&\hspace{-1.5mm} {\scriptstyle H_n(A_k)} \hspace{-1.5mm}&\hspace{-1.5mm} \stackrel{i_k}{\longrightarrow} \hspace{-1.5mm}&\hspace{-1.5mm} {\scriptstyle H_n(X_k)}
&\hspace{-1.5mm} \stackrel{j_k}{\longrightarrow} \hspace{-1.5mm}&\hspace{-1.5mm} {\scriptstyle H_n(X_k,A_k)} \hspace{-1.5mm}&\hspace{-1.5mm}\stackrel{\partial_k}{\longrightarrow} \hspace{-1.5mm}&\hspace{-1.5mm}
{\scriptstyle H_{n-1.5}(A_k)}\hspace{-1.5mm}&\hspace{-1.5mm} \longrightarrow \cdots\vspace{1mm}\\
&^{\alpha_k}\downarrow\quad&&^{\gamma_k}\downarrow\quad&&^{\beta_k}\downarrow\quad&&^{\alpha_k}\downarrow\quad\quad&&\\
\cdots\longrightarrow \hspace{-1.5mm}&\hspace{-1.5mm} {\scriptstyle H^{r_k-n}(X_k,A_k^c)} \hspace{-1.5mm}&\hspace{-1.5mm} \stackrel{q^*_k}{\longrightarrow} \hspace{-1.5mm}&\hspace{-1.5mm}{\scriptstyle H^{r_k-n}(X_k)}
\hspace{-1.5mm}&\hspace{-1.5mm} \stackrel{p^*_k}{\longrightarrow} \hspace{-1.5mm}&\hspace{-1.5mm}
{\scriptstyle H^{r_k-n}(A_k^c)} \hspace{-1.5mm}&\hspace{-1.5mm}\stackrel{\partial^*_k}{\longrightarrow} \hspace{-1.5mm}&\hspace{-1.5mm}
{\scriptstyle H^{r_k-n+1}(X_k,A^c_k)} \hspace{-1.5mm}&\hspace{-1.5mm}\longrightarrow\cdots
  \end{array}\]
we have $H^*(X_k,A_k^c)={\rm coim}\,q^*_k\oplus\Sigma^{-1}{\rm im}\,\partial^*_k
={\rm ker}\,p^*_k\oplus{\rm coker}\,p^*_k\subset
H^*_\XX(X_k,A_k^c)$.
So the following diagrams are commutative
\[\begin{array}{ccc}
H_*(A_k) &\subset&H_*^\XX(X_k,A_k)\,\vspace{1mm}\\
^{\alpha_k}\downarrow\quad\,&&^{\theta_k}\downarrow\quad\quad\\
H^*(X_k,A_k^c) &\subset&H^*_\XX(X_k,A_k^c),\end{array}\quad\quad
\begin{array}{ccc}
H_*(X_k) &\subset&H_*^\XX(X_k,A_k)\,\vspace{1mm}\\
^{\gamma_k}\downarrow\quad\,&&^{\theta_k}\downarrow\quad\quad\\
H^*(X_k) &\subset&H^*_\XX(X_k,A_k^c),\end{array}
\]
where $\theta_k,\alpha_k,\gamma_k$ are as in the proof of Theorem 5.1.
This implies that the following diagram is commutative
\[\begin{array}{ccc}
{\scriptstyle H_*(M)}
&\stackrel{\alpha_M}{-\!\!\!-\!\!\!-\!\!\!-\!\!\!-\!\!\!\longrightarrow}&
{\scriptstyle H^{r-*}(\w X,M^c)}\,\vspace{1mm}\\
\|&&\|\\
{\scriptstyle H_*(Y_1){\otimes}{\cdots}{\otimes}H_*(Y_m)}
&\stackrel{\alpha_M}{-\!\!\!-\!\!\!-\!\!\!-\!\!\!-\!\!\!\longrightarrow}&
{\scriptstyle H^{r_1-*}(X_1,Y_1^c){\otimes}{\cdots}{\otimes}H^{r_m-*}(X_m,Y_m^c)}\,\vspace{1mm}\\
\|\wr&&\|\wr\\
\oplus_{\sigma\subset S,\,\omega{\cap}S=\emptyset}\,{\scriptstyle H_*^{\sigma\!,\,\omega}(\underline{X},\underline{A})}
&\stackrel{\oplus\,\Sigma^{|\omega|}\gamma_{\sigma\!,\omega}}{-\!\!\!-\!\!\!-\!\!\!-\!\!\!-\!\!\!\longrightarrow}&
\oplus_{\tilde\sigma\subset S,\,\omega{\cap}S=\emptyset}\,{\scriptstyle
\Sigma H^{r-*}_{\tilde \sigma\!,\,\omega}(\underline{X},\underline{A}^c)}\vspace{1mm}\\
\cap&&\cap\\
{\scriptstyle H_*^\XX(X_1,A_1){\otimes}{\cdots}{\otimes}H_*^\XX(X_m,A_m)}
&\stackrel{\theta_1\otimes\cdots\otimes\theta_m}
{-\!\!\!-\!\!\!-\!\!\!-\!\!\!-\!\!\!\longrightarrow}&
{\scriptstyle \Sigma(H^*_\XX(X_1,A_1^c){\otimes}{\cdots}{\otimes}H^*_\XX(X_m,A_m^c))},
\end{array}\quad(4)\]
where the $\Sigma^{|\omega|}$ of $\gamma_{\sigma\!,\,\omega}$ comes from
the desuspension isomorphism $\Sigma^{-1}{\rm im}\,\partial^*_k\cong{\rm coker}\,p^*_k$
and the $\Sigma$ of $H^*_\XX(X_1,A_1^c){\otimes}{\cdots}{\otimes}H^*_\XX(X_m,A_m^c)$ comes from
the isomorphism $\overline H^*(\w X,M^c)\cong\Sigma\overline H^*(M^c)$.

For $\sigma\subset S$, $H_*^{\sigma\!,\,\omega}(\Delta\!^S)=0$ if $\omega{\cap}S\neq\emptyset$ and
$H_*^{\sigma\!,\,\omega}(\Delta\!^S)\cong{\Bbb Z}$ if $\omega{\cap}S=\emptyset$.
So $\gamma_{\Delta\!^S\!,\sigma\!,\,\omega}\colon 0\to 0$ if $\omega{\cap}S\neq\emptyset$
and $\gamma_{\Delta\!^S\!,\sigma\!,\,\omega}\colon\Bbb Z\to\Bbb Z$ if $\omega{\cap}S=\emptyset$.
For $\omega{\cap}S=\emptyset$, identify $H_*^{\sigma\!,\,\omega}(K){\otimes}H_*^{\sigma\!,\,\omega}(\underline{X},\underline{A})$ and
$H^*_{\tilde\sigma\!,\,\omega}(K^\circ){\otimes}H^*_{\tilde\sigma\!,\,\omega}(\underline{X},\underline{A}^c)$ respectively with
$H_*^{\sigma\!,\,\omega}(\underline{X},\underline{A})$ and $\Sigma^{|\omega|-1}\,H^*_{\tilde\sigma\!,\,\omega}(\underline{X},\underline{A}^c)$, then we
have the following commutative diagram
\[\begin{array}{ccc}
{\scriptstyle H_0^{\sigma\!,\,\omega}(\Delta\!^S){\otimes}H_*^{\sigma\!,\,\omega}(\underline{X},\underline{A})}
&\stackrel{\gamma_{\Delta\!^S\!,\sigma,\omega}\otimes\gamma_{\sigma,\omega}}
{-\!\!\!-\!\!\!-\!\!\!-\!\!\!-\!\!\!\longrightarrow}&
{\scriptstyle H^{|\omega|-1}_{\tilde\sigma\!,\,\omega}((\Delta\!^S)^\circ){\otimes}
H^{r-|\omega|-*}_{\tilde\sigma\!,\,\omega}(\underline{X},\underline{A}^c)}
\,\vspace{1mm}\\
\|\wr&&\|\wr\\
{\scriptstyle H_*^{\sigma\!,\,\omega}(\underline{X},\underline{A})}
&\stackrel{\Sigma^{|\omega|-1}\gamma_{\sigma\!,\,\omega}}{-\!\!\!-\!\!\!-\!\!\!-\!\!\!-\!\!\!\longrightarrow}&
{\scriptstyle
H^{r-*-1}_{\tilde \sigma\!,\,\omega}(\underline{X},\underline{A}^c)}.
\end{array}\]

The direct sum of the above isomorphisms for all $\sigma\subset S$ and $\omega{\cap}S=\emptyset$ is just
the third row of (4).
$\hat\alpha_M$ is the direct sum of
the above isomorphisms for all $\sigma\subset S$ and $\omega=\emptyset$.
$\overline\alpha_M$ is the direct sum of the above isomorphisms for all $\sigma\subset S$, $\omega\neq\emptyset$ and
$\omega{\cap}S=\emptyset$.
So $\alpha_M=\hat\alpha_M{\oplus}\overline\alpha_M$ for the special case $M={\cal Z}(\Delta\!^S;\underline{X},\underline{A})$.
\hfill $\Box$\vspace{3mm}

{\bf Example 5.7}\, Regard $S^{r+1}$ as one-point compactification of $\Bbb R^{r+1}$.
Then for $q\leqslant r$, the standard  space pair $(S^{r+1},S^q)$ is given by\vspace{1mm}\\
\hspace*{2mm}$S^q=\{(x_1,{\cdots},x_{r+1})\in\Bbb R^{r+1}\subset S^{r+1}
\,|\,x^2_1{+}{\cdots}{+}x_{q+1}^2=1,\,x_i=0,\,\,{\rm if}\,\,i>q{+}1\}.$

Let $M={\cal Z}_K\Big(\!
\begin{array}{ccc}
\scriptstyle{r_1{+}1} \!&\!\scriptstyle{\cdots}\!&\!\scriptstyle{r_m{+}1}\\
\scriptstyle{ q_1}   \!&\!\scriptstyle{\cdots}\!&\! \scriptstyle{q_m}\end{array}\!\Big)
={\cal Z}(K;\underline{X},\underline{A})$ be the polyhedral product space
such that $(X_k,A_k)=(S^{r_k+1},S^{q_k})$.
Since  $S^{r-q}$ is a deformation retract of $S^{r+1}{\setminus}S^q$,
the complement space $M^c={\cal Z}(K^\circ;\underline{X},\underline{A}^c)$
is homotopic equivalent to ${\cal Z}_{K^\circ}\Big(
\begin{array}{ccc}
\scriptstyle{r_1+\,1\,} &\scriptstyle{\cdots}&\scriptstyle{r_m+\,1\,}\\
\scriptstyle{r_1{-}q_1}   &\scriptstyle{\cdots}& \scriptstyle{r_m{-}q_m}\end{array}\Big)$.\vspace{1mm}

Since all $H_*^{\sigma\!,\,\omega}(\underline{X},\underline{A})\cong {\Bbb Z}$,
we may identify $H_*^{\sigma\!,\,\omega}(K){\otimes}H_*^{\sigma\!,\,\omega}(\underline{X},\underline{A})$
with $\Sigma^t H_*^{\sigma\!,\,\omega}(K)$,
where $t=\Sigma_{k\in\sigma}(r_k{+}1)+\Sigma_{k\in\omega}\,q_k$.
For $\sigma\subset[m]$, let ${\Bbb Z}_\sigma$ be the free group generated by
$\sigma$ with degree $0$. Then
\[\hat H_*(M)=\oplus_{\sigma\in K}\,\Sigma^{\Sigma_{k\in\sigma}\,(r_k{+}1)}{\Bbb Z}_\sigma,\]
\[\overline H_*(M)=\oplus_{(\sigma\!,\,\omega)\in\overline \XX_m}\,
\Sigma^{\Sigma_{k\in\sigma}\,(r_k{+}1)+\Sigma_{k\in\omega}\,q_k}
H_*^{\sigma\!,\,\omega}(K).\]

Dually, the cohomology of the complement space $M^c$ is
\[\hat H^*(M^c)=\oplus_{\sigma\in K^\circ}\,\Sigma^{\Sigma_{k\in\sigma}\,(r_k{+}1)}{\Bbb Z}_\sigma,\]
\[\overline H^*(M^c)=\oplus_{(\sigma\!,\,\omega)\in\overline \XX_m}\,
\Sigma^{\Sigma_{k\in\sigma}\,(r_k{+}1)+\Sigma_{k\in\omega}\,(r_k-q_k)}
H^*_{\sigma\!,\,\omega}(K^\circ).\]

In this case, the direct sum of $\gamma_{K,\sigma\!,\omega}
\colon H_*^{\sigma\!,\,\omega}(K)\to H^{|\omega|-*-1}_{\tilde\sigma\!,\,\omega}(K^\circ)$
over all $(\sigma\!,\omega)\in\overline\XX_m$ (regardless of degree) is the isomorphism
$\overline H_*(M)\cong\overline H\,^{r-*-1}(M^c)$.

Specifically, ${\cal Z}(K;S^{2n+1},S^n)={\cal Z}_K\Big(\!
\begin{array}{ccc}
\scriptstyle{2n{+}1} \hspace{-1mm}&\hspace{-1mm}\scriptstyle{\cdots}\hspace{-1mm}&\hspace{-1mm}\scriptstyle{2n{+}1}\\
\scriptstyle{n}   \hspace{-1mm}&\hspace{-1mm}\scriptstyle{\cdots}\hspace{-1mm}&\hspace{-1mm} \scriptstyle{n}\end{array}\!\Big)$.
Then we have $$\overline H_*({\cal Z}(K;S^{2n+1},S^n))\cong
\overline H\,^{(2n+1)m-*-1}({\cal Z}(K^\circ;S^{2n+1},S^n)).$$

\end{document}